\documentclass[final,onefignum,onetabnum]{siamonline181217}

\usepackage{tikz}
\usepackage[shortlabels]{enumitem}
\usepackage{pgfplots}
\usepackage[algo2e]{algorithm2e} 

\usetikzlibrary{arrows}
\usepackage{float}
\graphicspath{{RadonDensityFigs2/}}
\usetikzlibrary{decorations.pathmorphing}
\tikzset{snake it/.style={decorate, decoration=snake}}
\newcount\gaussF
\edef\gaussR{0}
\edef\gaussA{0}

\makeatletter
\pgfmathdeclarefunction{gaussR}{0}{%
 \global\advance\gaussF by 1\relax
 \ifodd\gaussF
  \pgfmathrnd@%
  \ifdim\pgfmathresult pt=0.0pt\relax%
    \def\pgfmathresult{0.00001}%
  \fi
  \pgfmathln@{\pgfmathresult}%
  \pgfmathmultiply@{-2}{\pgfmathresult}%
  \pgfmathsqrt@{\pgfmathresult}%
  \global\let\gaussR=\pgfmathresult
  \pgfmathrnd@%
  \pgfmathmultiply@{360}{\pgfmathresult}%
  \global\let\gaussA=\pgfmathresult
  \pgfmathcos@{\pgfmathresult}%
  \pgfmathmultiply@{\pgfmathresult}{\gaussR}%
 \else
  \pgfmathsin@{\gaussA}%
  \pgfmathmultiply@{\gaussR}{\pgfmathresult}%
\fi
}

\pgfmathdeclarefunction{invgauss}{2}{%
  \pgfmathln{#1}
  \pgfmathmultiply@{\pgfmathresult}{-2}%
  \pgfmathsqrt@{\pgfmathresult}%
  \let\@radius=\pgfmathresult%
  \pgfmathmultiply{6.28318531}{#2}
  \pgfmathdeg@{\pgfmathresult}%
  \pgfmathcos@{\pgfmathresult}%
  \pgfmathmultiply@{\pgfmathresult}{\@radius}%
}

\pgfmathdeclarefunction{randnormal}{0}{%
  \pgfmathrnd@
  \ifdim\pgfmathresult pt=0.0pt\relax%
    \def\pgfmathresult{0.00001}%
  \fi%
  \let\@tmp=\pgfmathresult%
  \pgfmathrnd@%
  \ifdim\pgfmathresult pt=0.0pt\relax%
    \def\pgfmathresult{0.00001}%
  \fi
  \pgfmathinvgauss@{\pgfmathresult}{\@tmp}%
}

\DeclareMathOperator*{\argmin}{arg\,min}
\DeclareMathOperator{\sinc}{sinc}

\usepackage{lipsum}
\usepackage{amsfonts}
\usepackage{graphicx}
\usepackage{subcaption}
\usepackage{epstopdf}
\usepackage{algorithmic}
\ifpdf
  \DeclareGraphicsExtensions{.eps,.pdf,.png,.jpg}
\else
  \DeclareGraphicsExtensions{.eps}
\fi

\usepackage{enumitem}
\setlist[enumerate]{leftmargin=.5in}
\setlist[itemize]{leftmargin=.5in}

\newsiamremark{remark}{Remark}
\newsiamremark{hypothesis}{Hypothesis}
\crefname{hypothesis}{Hypothesis}{Hypotheses}
\newsiamthm{claim}{Claim}

\headers{Non-parametric density estimation}{J. Webber, E. Hussey, E. Miller and S. Aeron}

\title{On non-parametric density estimation on linear and non-linear manifolds using generalized Radon transforms\thanks{Submitted to the editors DATE.
\funding{This work was funded by Center for Applied Brain and Cognitive Sciences (CABCS) at Tufts University. Shuchin Aeron was supported in part by the NSF CAREER award.}}}

\author{James Webber\thanks{Dept. of ECE, Tufts University 
  (\email{James.Webber@tufts.edu})}
\and Erika Hussey \thanks{Center for Applied Brain and Cognitive Sciences (CABCS)
  (\email{erika.hussey@tufts.edu})}
\and Eric Miller \thanks{Dept. of ECE, Tufts University
  (\email{elmiller@ece.tufts.edu})}
\and Shuchin Aeron\thanks{Dept. of ECE, Tufts University
  (\email{shuchin@ece.tufts.edu})}
  }

\usepackage{amsopn}

\begin{document}

\maketitle

\begin{abstract}
Here we present a new non-parametric approach to density estimation and classification derived from theory in Radon transforms and image reconstruction. We start by constructing a ``forward problem" in which the unknown density is mapped to a set of one dimensional empirical distribution functions computed from the raw input data.  Interpreting this mapping in terms of Radon-type projections provides an analytical connection between the data and the density with many very useful properties including stable invertibility, fast computation, and significant theoretical grounding. Using results from the literature in geometric inverse problems we give uniqueness results and stability estimates for our methods. We subsequently extend the ideas to address problems in manifold learning and density estimation on manifolds. We introduce two new algorithms which can be readily applied to implement density estimation using Radon transforms in low dimensions or on low dimensional manifolds embedded in $\mathbb{R}^d$. We test our algorithms performance on a range of synthetic 2-D density estimation problems, designed with a mixture of sharp edges and smooth features. We show that our algorithm can offer a consistently competitive performance when compared to the state--of--the--art density estimation methods from the literature.
\end{abstract}



\section{Introduction}

\begin{figure}
\centering
\begin{subfigure}{0.42\textwidth}
\includegraphics[width=1.0\linewidth, height=5.5cm]{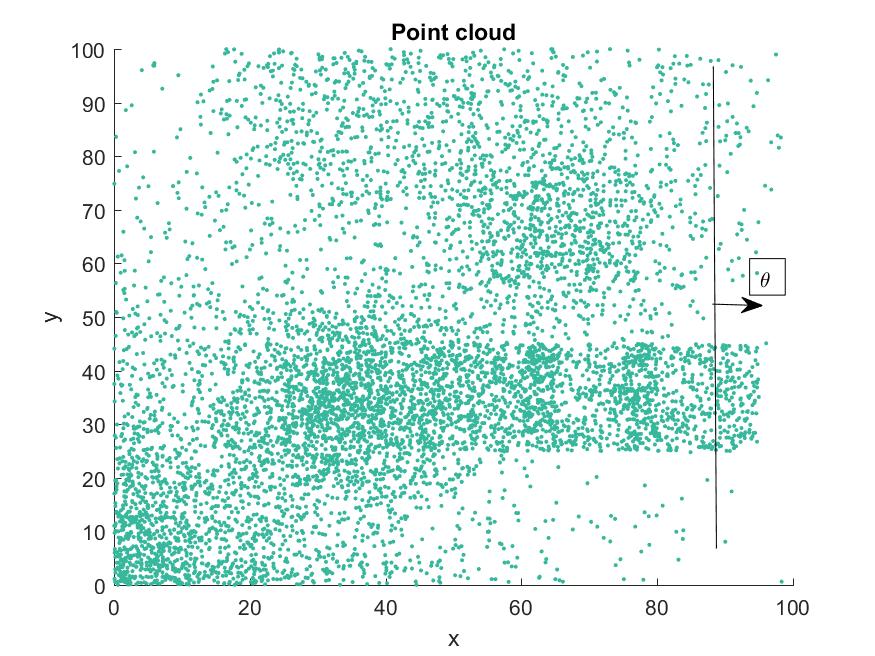} 
\end{subfigure}\hspace{5mm}
\begin{subfigure}{0.42\textwidth}
\includegraphics[width=1.0\linewidth, height=5.5cm]{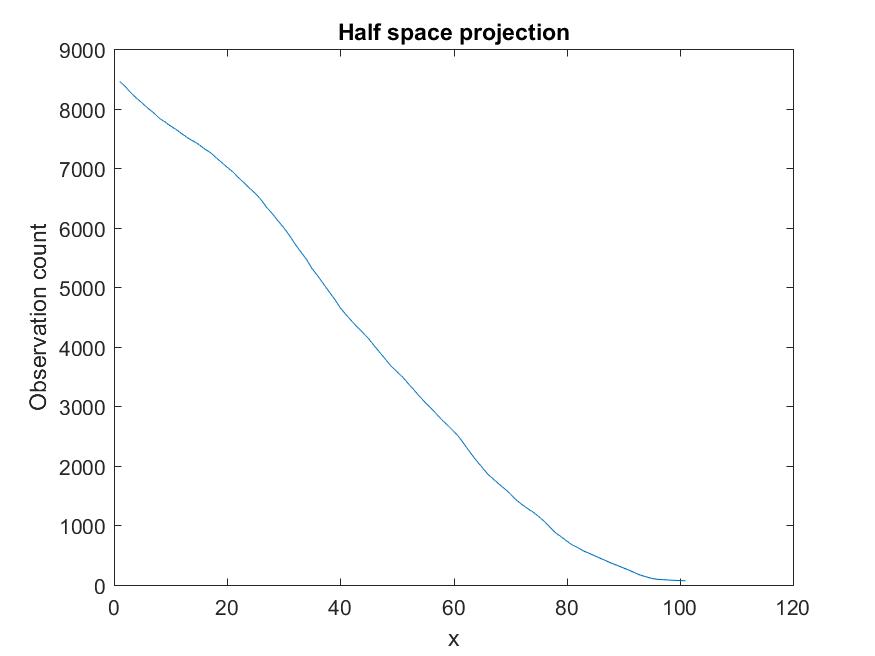} 
\end{subfigure}
\caption{Point cloud of IID samples from an unknown density $f$ (left). Half space projection of $f$ along the $x$ axis (right). Here $\theta=(1,0)$ is in the direction of the $x$ axis. The values in the right hand plot are the density count in the right half space to the line shown in the left hand plot as it is translated in the direction $\theta$.}
\label{fig1.1}
\end{figure}

\textbf{Motivation}: Consider Figure \ref{fig1.1}, a point cloud that is obtained by IID sampling from an unknown density $f : \mathbb{R}^n\to\mathbb{R}$ (in this case $n=2$, but we consider the theoretical aspects of the general $n$ dimensional case in this paper). A principal goal in machine learning is to estimate the distribution function $f$ from a given finite sample, under the assumption that $f$ is smooth \cite{kerselec,gdense,BI,TV1D,TV1D1}. The primary contribution of this paper is the development and analysis of methods for recovering $f$ from such a sample set using ideas from the theory of Radon transforms.  More specifically, Radon theory provides for a highly flexibly  family of algorithms for stably estimating $f$ from projections of the empirical data set onto low dimensional manifolds such as affine subspaces, spheres, and the like. See figure \ref{fig1.1} for an approximation of a half space Radon transform projection using the empirical distribution function (by observation counting in half spaces). Moreover, unlike more common density estimation methods, the Radon-based technique is naturally adapted to densities, which may not be smooth. The prediction of one dimensional densities is well studied, from a theoretical standpoint \cite{DKW} and in implementation. Indeed there are many effective techniques which can be applied with good results in the one dimensional case \cite{gdense,kerselec,TV1D,TV1D1}. It is well known that the problem of reconstructing a density $f$ (of a certain type) from its Radon transform $Rf$, given enough projection data, is mildly ill--posed \cite[page 42]{natterer}, and hence we can expect a small amplification in the error from a prediction of $Rf$ in the reconstruction of $f$ from $Rf$. By combining known approximation techniques for one dimensional density estimation and reconstruction methods in image reconstruction and inverse problems, we propose new techniques in non-parametric density estimation, for which we can derive convergence estimates and demonstrate the effectiveness of our methods in implemenation on a variety of challenging synthetic density estimation problems. 

\section{Relevant literature}
The classical and still the most-widely used methods for non-parametric density estimation are based on a Kernel Density Estimate (KDE) \cite{kerselec,gdense,BI,Fukunaga:1975up}, which approximates the unknown density by its convolution with a known density function (e.g. isotropic Gaussians). Other methods in non-parametric density estimation consider the reconstruction of the density by maximizing the log--likelihood $\max_f \sum_i \log f(X_i)$ \cite{TV1D,TV1D1} with additional smoothness regularization on the density $f$. We give a detailed comparison of such ideas and our own approach later in section \ref{method}. Extensions of the KDE approach to density estimation on manifolds (embedded in high dimensions) has been considered in \cite{Kuleshov:2017vp,Pelletier:2005fk,Ozakin:2009wq}.

The approaches based on KDE can fail to capture precise details in the density (e.g. edge effects) due to smoothing, particularly if the bandwidth is too large. To alleviate this issue the use of the hyperplane Radon transform for density estimation has been considered in the literature \cite{rad1,rad2,rad3}. In \cite{rad2} the hyperplane Radon transform is approximated by the application of a smoother to empirical density functions and the filtered backprojection formula is used to reconstruct Gaussian mixture densities in two dimensions. In \cite{rad1} it is shown that when the Radon transform is approximated by a series of one dimensional kernels, the reconstruction is a type of convolution kernel estimate. Using this idea new kernel density estimators are derived. In \cite{rad3} a new approach to Gaussian mixture model fitting is presented which uses one dimensional projections of the target density, which are described by the Radon transform, to give a more stable convergence with random initializations than the Expectation Maximization (EM) algorithm.

\section{Main contributions} This paper advances the use of Radon transforms for non-parametric density estimation in the following novel directions:
\begin{itemize}
\setlength \itemsep{2pt}
    \item We use new types of empirical distribution functions to approximate generalized Radon transform projections of the density. Here we consider spherical and half space Radon transforms, whereas hyperplane transforms are only considered in previous work. This allows us to elegantly derive error estimates for our approach and provides the basis for a new, fast method in density estimation, which is completely non-parametric.
    \item While the current methods use the explicit inversion formula for the hyperplane Radon transform (a filtered backprojection approach), we discretize generalized Radon transform operators on pixel grids (assuming a piecewise constant density) and use existing machinary for solving regularized inverse problems with Total Variation (TV) regularization. We choose TV as it is a powerful regularizer commonly applied image reconstruction which favours piecewise smooth solutions \cite{TV1,TV2}. The idea of using TV for density estimation has so far been only generally applied to the 1-D case when using the likelihood approach \cite{TV1D,TV1D1}. A 2-D example is considered in \cite{TV1D1} but we extend on the bivariate case significantly and introduce a variety of synthetic 2-D test examples in this paper.
\item For the test densities considered, we give a full error and image quality comparison using spherical Radon transforms, half space Radon transforms, the methods presented in \cite{rad2,TV1D1} and a kernel estimate. We see a consistently high performance using spherical Radon transforms when compared to the other methods considered, and across a range of point cloud sample sizes.
    \item We derive convergence estimates for our approach by combining the theory of Radon transforms \cite[page 94]{natterer} and results on the estimation of one dimensional densities by empirical distributions as described by the Dvoretzky--Kiefer--Wolfowitz inequality \cite{DKW}. In \cite{rad2} convergence estimates are given for a hyperplane approximation if certain assumptions are met. For assumption (c) in 5.1.2. in \cite{rad2} to hold based on the suggested approximation using truncated Fourier series, we require knowledge of the Sobolev norms $\|g_{\theta}\|_k$, where $g_{\theta}$ are the (unknown) projections to be approximated. Here we give estimates for half space and spherical approximators only assuming a degree of smoothness for the underlying distribution.
    \item Somewhat in a vein similar to extending KDE methods to density estimation on manifolds,  \cite{Kuleshov:2017vp,Pelletier:2005fk,Ozakin:2009wq}, we extend the Radon idea from Euclidean space to density estimation on manifolds embedded in $\mathbb{R}^d$ by reconstructing the density on local tangent spaces. This combines the theory of Radon transforms with the theory of manifold sampling and reconstruction, taking inspiration from a number of methods in non linear dimensionality reduction, such as Locally Linear Embedding \cite{lle} and Local Tangent Space Alignment (LTSA) \cite{LTSA,fefferman}. The main idea is to embed the data locally into low-dimensions using these methods and then use the proposed density estimation approach on these patches.
\item We present two new algorithms which can be readily applied to implement our density estimation techniques in low dimensional Euclidean space and on low dimensional manifolds embedded in $\mathbb{R}^d$. The codes are available from the authors on request.
\end{itemize}

\section{Organisation of the paper} This paper is organized as follows. In section \ref{mainsec} we state some mathematical preliminaries and theory on Radon transform, and explain the intuitive microlocal idea behind edge detection and reconstruction using Radon transforms. We then introduce our main idea, where we propose to use empirical distribution functions derived from Independent and Identically Distributed (IID) point clouds to approximate Radon transforms of the probability density.

In section \ref{error1} we provide error and convergence estimates for our approach by combining theory from geometric inverse problems \cite[page 94]{natterer} and empirical distribution function estimators \cite{DKW}. We show that the convergence rate of the least squares error in our approximate density is $(1/m)^{\frac{1}{2(n+2)}}$ when using the half space Radon transform, where $m$ is the number of observations and $n$ is the dimension. We also give estimates for the expected error in a spherical Radon transform approximation and explain the increase in error with $n$, the number of projections taken and the confidence level.

In section \ref{method} we present simulated density reconstructions in two dimensions using an algebraic approach with TV regularization, and give an error and side by side image reconstruction comparison with the methods presented in \cite{rad2,TV1D1} and a Gaussian kernel estimate. The algebraic approach implemented is formalized in Algorithm \ref{alg1} for density estimation by inverse Radon transform in low dimensional Euclidean space. Later in section \ref{manifolds} we extend our ideas from Euclidean space in low dimensions to density estimation on low dimensional manifolds embedded in $\mathbb{R}^d$, using theory and techniques in manifold learning \cite{fefferman,isomap,lle}. 


\section{The main idea and approach}
\label{mainsec}
In this section we state some mathematical preliminaries and theory on Radon transforms, and explain briefly their inversion and edge detection properties. We then introduce our main idea for using empirical distribution functions of IID point clouds to approximate Radon transforms of the probability density. 
\subsection{Mathematical preliminaries and Radon Transforms}
Throughout this paper $X=\{x_1,\ldots,x_m\}$ will denote a point could in $\mathbb{R}^n$ and $f$ will denote the density from which $X$ is drawn. We also denote:
\begin{enumerate}
\setlength \itemsep{3pt}
\item $Z=\mathbb{R}\times S^{n-1}$ as the unit cylinder in $\mathbb{R}^{n+1}$.
\item $C^{k}_0(\Omega)$ as the space of $k$ continuously differentiable functions compactly supported on $\Omega\subset\mathbb{R}^n$.
\item $\Omega_n$ as the unit ball in $\mathbb{R}^n$.
\item $\|g\|^p_{L^p(\mathbb{R})}=\int_{S^{n-1}}\int_{\mathbb{R}}|g(s,\theta)|^p\mathrm{d}s\mathrm{d}\theta$ as the $L^p$ norm for functions $g$ on the cylinder.
\item $V_{s,n-1}$ and $w_{s,n-1}$ as the volume and surface area of a sphere radius $s$ in $\mathbb{R}^n$ respectively.
\end{enumerate}
\begin{definition}[{\cite[page 9]{natterer}}]
Let $f\in C^{\infty}_0(\mathbb{R}^n)$. Then the hyperplane Radon transform of $f$ is defined as
\begin{equation}
\label{equ7.1}
Rf(s,\theta)=\int_{\mathbb{R}^n}f(x)\delta(x\cdot\theta-s)\mathrm{d}x,
\end{equation}
where $\delta$ denotes the Dirac delta function. 
\end{definition}
\begin{definition}
Let $f\in C^{\infty}_0(\mathbb{R}^n)$. Then the spherical (volumetric) Radon transform of $f$ is defined as
\begin{equation}
\label{equ6.1}
R_Mf(s,x)=\int_{0}^{s}r^{n-1}\int_{S^{n-1}}f(x+r\theta)\mathrm{d}\theta\mathrm{d}r.
\end{equation}
\end{definition}
Here the set of spheres in $\mathbb{R}^n$ are parameterized by $(s,x)$, where $s\in\mathbb{R}^+$ is the radius and $x\in\mathbb{R}^n$ is the center of the sphere. After differentiating with respect to the $s$ variable and dividing by $s^{n-1}$ the transform $R_Mf$ is equivalent to the spherical means transform defined in \cite{quinto}.
\begin{definition}
Let $f\in C^{\infty}_0(\mathbb{R}^n)$. Then the half space Radon transform of $f$ is defined as
\begin{align}
\label{equ7.11}
R_Hf(s,\theta)=\int_{-\infty}^{s}\int_{\mathbb{R}^n}f(x)\delta(x\cdot\theta-r)\mathrm{d}x\mathrm{d}r,
\end{align}
where $\delta$ denotes the Dirac delta function. 
\end{definition}

\begin{figure}[!h]
\centering
\begin{tikzpicture}[scale=3]
\draw [->] (0,-0.5)--(0,1);
\draw [->] (-1,0)--(1,0);
\draw [->] (0,0)--(0.25,0.25);
\draw (-0.1,0.6)--(0.6,-0.1);
\draw (-0.5,0.5) circle (0.3);
\node at (0.3,0.3) {$\theta$};
\node at (0.12,0.18) {$s$};
\node at (-0.55,0.5) {$x$};
\draw [<->] (-0.5,0.5)--(-0.5,0.8);
\node at (-0.55,0.65) {$s$};
\draw[fill] (-0.5,0.5) circle (0.01);
\end{tikzpicture}
\caption{A hyperplane in two dimensions parameterized by the perpendicular distance $s$ from the origin and a direction $\theta$ (right). A circle parameterized by a center point $x$ and radius $s$ (left).}
\label{fig2}
\end{figure}
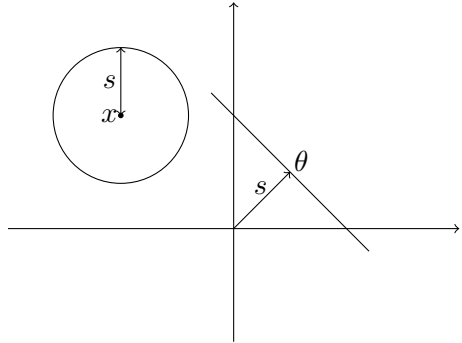

Here the set of hyperplanes in $\mathbb{R}^n$ are parameterized by $(s,\theta)\in Z$. See figure \ref{fig2}. After differentiating with respect to the $s$ variable the transform $R_Hf$ reduces to the classical hyperplane Radon transform $Rf$.

It is well known \cite[page11, Theorem 2.5]{natterer1}, that the hyperplane (and hence half space) Radon transforms have an explicit left inverse and are injective for $f\in C^{\infty}_0(\mathbb{R}^n)$. Similarly for the spherical Radon transform we can derive explicit reconstruction formulae for the density \cite{quinto,halt1,halt2}, and the solution is unique if enough spherical centers and radii are known (e.g. for any fixed radius $s$, if all centers $x$ are known, the reconstruction is obtained by spherical deconvolution or deblurring). The previous works of \cite{rad2} apply the formula \cite[page11, Theorem 2.5]{natterer1} directly in density estimation. Our reconstruction method discretizes the linear operator $R_H$ (or $R_M$) on voxel grids and solves the resulting set of sparse linear equations with regularization.

\begin{figure}[!h]
\centering
\begin{subfigure}{0.35\textwidth}
\includegraphics[width=1\linewidth, height=4.8cm]{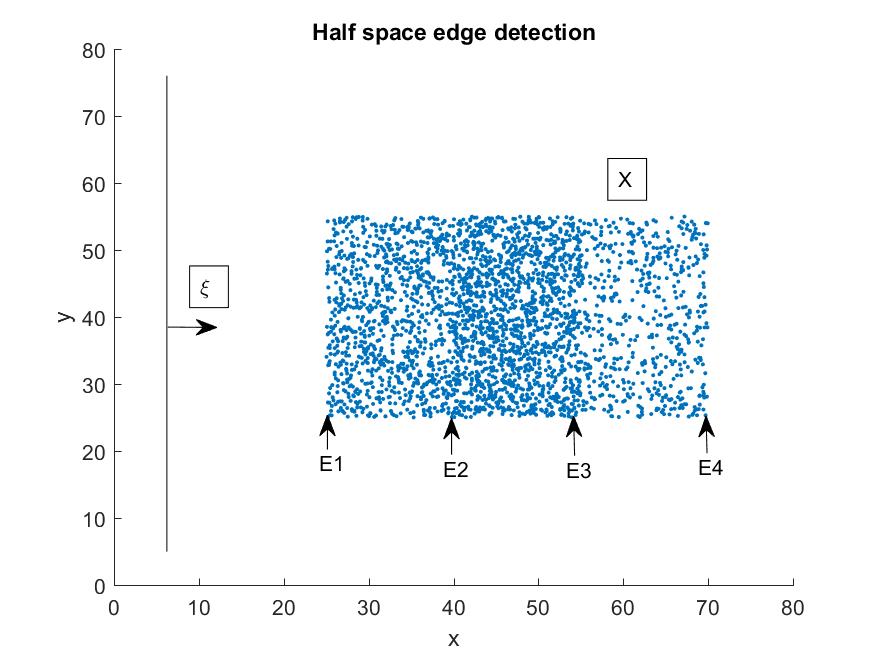} 
\end{subfigure}
\begin{subfigure}{0.35\textwidth}
\includegraphics[width=1\linewidth, height=4.8cm]{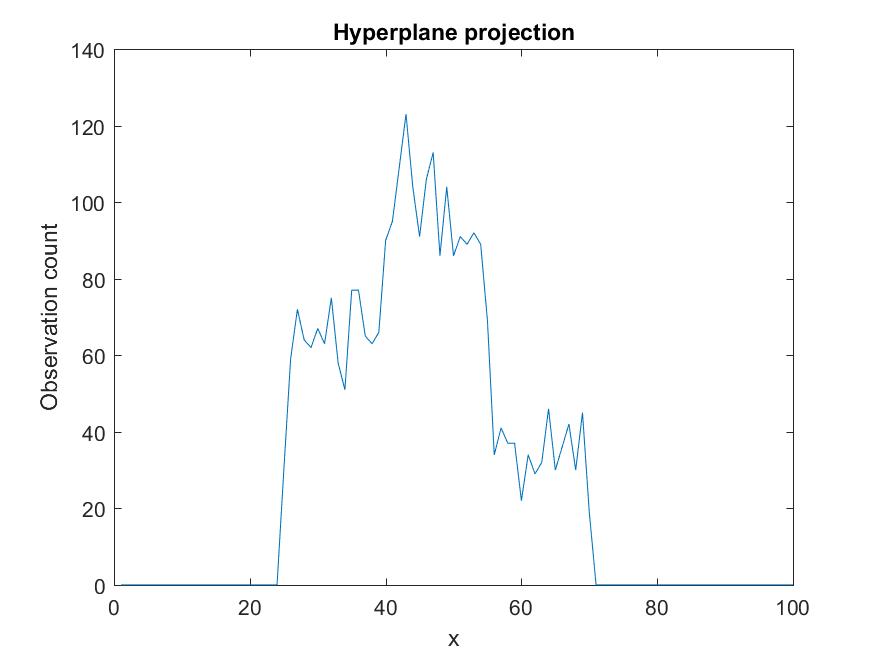}
\end{subfigure}
\begin{subfigure}{0.35\textwidth}
\includegraphics[width=1\linewidth, height=4.8cm]{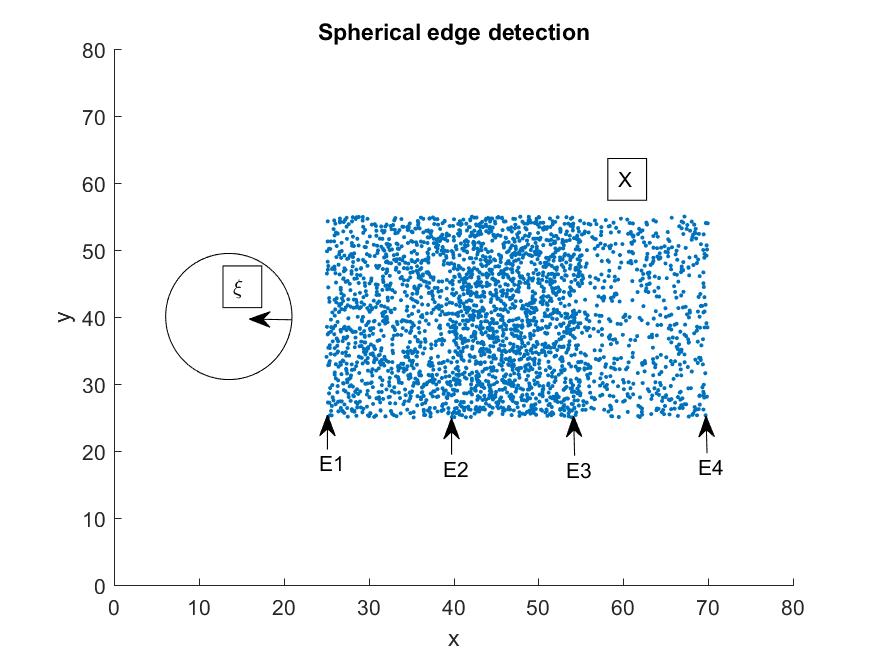} 
\end{subfigure}
\begin{subfigure}{0.35\textwidth}
\includegraphics[width=1\linewidth, height=4.8cm]{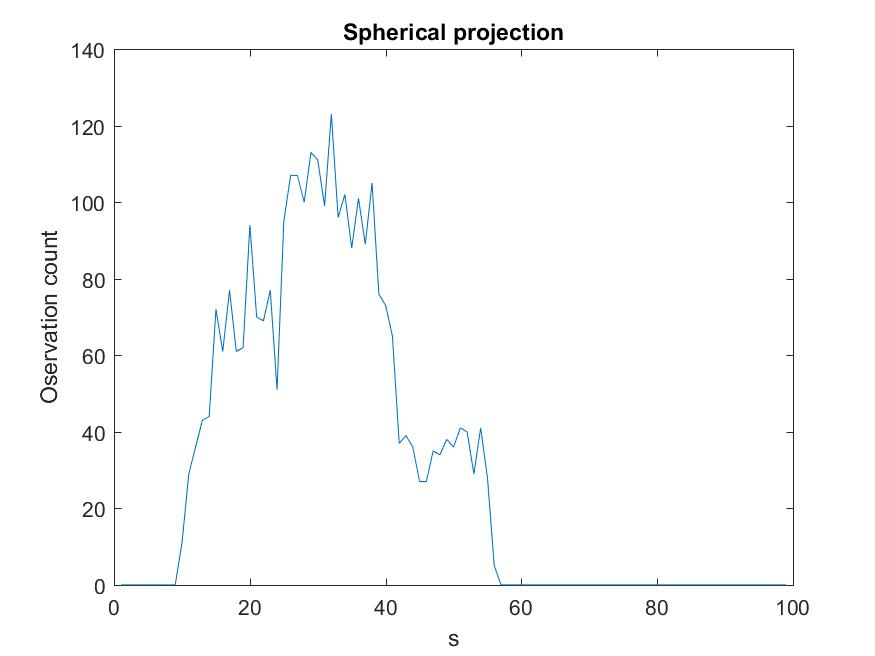}
\end{subfigure}
\caption{An IID sample taken from two overlapping uniform densities (left figures). There are four edges (marked by E1, E2, E3 and E4) in the direction $\xi$ highlighted. The jump discontinuities are present in the derivative of a spherical (bottom right) or half space (top right) projection. Here $s$ is the radius of the sphere with fixed center (as pictured). That is, with center $(x,y)=(12.5,40)$.}
\label{fig9.1}
\end{figure}

Radon transforms are known to have edge detection properties \cite{quinto2,quinto3}. If we know the integrals of $f$ over a set of hypersurfaces tangent to a direction $\xi\in S^{n-1}$ in a small neighborhood of a given $x\in\mathbb{R}^n$, then we can detect jump singularities in $f$ near $x$ in the direction $\xi$ (this forms the intuition behind the theory of microlocal analysis applied to Radon transforms). See figure \ref{fig9.1}. We can see that the jump discontinuities in the density are present in the derivatives (finite differences) of Radon transform (either spherical or half space) approximations, but the projections are noisy. We will see later that we can combat this noise while effectively preserving edges in the density by a discrete inversion of the spherical or half space Radon transform with a Total Variation (TV) penalty. TV is a powerful regularizer commonly applied in image reconstruction which favours piecewise smooth solutions.
\subsection{Application to density estimation: The core idea}
\begin{figure}[!h]
\centering
\begin{tikzpicture}[scale=0.8,
    declare function={a(\x)=0.75*\x-1.5;},
    declare function={b(\x)=0.75*\x-1;},
    declare function={c(\x)=4;},
    declare function={d(\x)=-4;},
    declare function={f(\x)=0.75*\x-1.25;},
]
\begin{axis}[
    domain=-4:4,
    axis lines=middle,
    axis equal image,
    xtick=\empty, ytick=\empty,
    enlargelimits=true,
    clip mode=individual, clip=false
]
\addplot [red, only marks, mark=*, samples=100, mark size=0.75]
    {0.5*(c(x)+d(x)) + 0.5*rand*(c(x)-d(x))};
\addplot [thick,domain=-4:4.5] {a(x)};
\coordinate (center) at (axis cs:-4.1,3);
\coordinate (c1) at (axis cs:6,3.3);
\coordinate (c2) at (axis cs:5.7,2.6);
\coordinate (c3) at (axis cs:5.27,3.15);
\end{axis}
\node at (center) {$X\subset\Omega_n$};
\end{tikzpicture}
\hspace{1cm}
\begin{tikzpicture}[scale=0.8,
    declare function={a(\x)=0.75*\x-1.5;},
    declare function={b(\x)=0.75*\x-1;},
    declare function={c(\x)=4;},
    declare function={d(\x)=-4;}
]
\begin{axis}[
    domain=-4:4,
    axis lines=middle,
    axis equal image,
    xtick=\empty, ytick=\empty,
    enlargelimits=true,
    clip mode=individual, clip=false
]
\addplot [red, only marks, mark=*, samples=100, mark size=0.75]
    {0.5*(c(x)+d(x)) + 0.5*rand*(c(x)-d(x))};
\coordinate (center1) at (axis cs:1,1);
\coordinate (center2) at (axis cs:1,1);
\coordinate (center) at (axis cs:4,5.6);
\coordinate (center4) at (axis cs:4.2,5.8);
\coordinate (center3) at (axis cs:1,3.1);
\end{axis}
\draw[thick] (center1) circle (1cm);
\end{tikzpicture}
\caption{Half space and spherical Radon transforms in two dimensions, approximated by counting observations in a set of half spaces and (inside of) spheres in the plane.}
\label{fig1}
\end{figure}
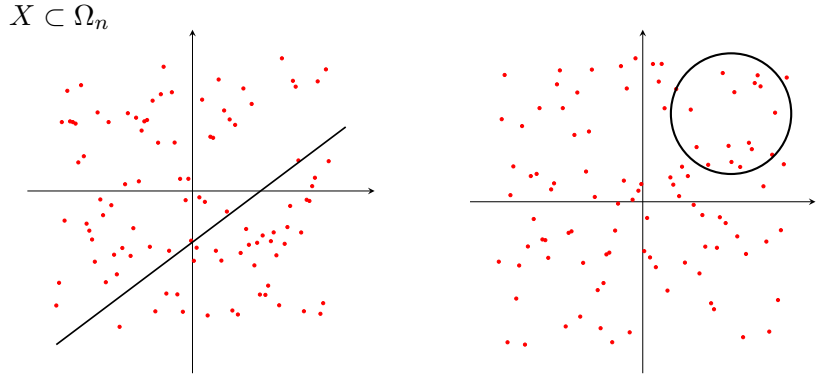

Let  $X=\{x_1,\ldots,x_m\}$ be a point cloud in $\Omega_n$ drawn as IID samples from a density $f\in C^{\infty}_0(\Omega_n)$. Then the half space Radon transform $R_H$ of $f$ can be approximated as the set of one dimensional empirical distribution functions
\begin{equation}
\begin{split}
\label{equ2.2}
R_Hf(s,\theta)=\int_{-\infty}^{s}\int_{\mathbb{R}^n}f(x)\delta(x\cdot\theta-r)\mathrm{d}x\mathrm{d}r\approx\sum_{i\in\{x_i\cdot\theta\leq s\}}1.
\end{split}
\end{equation}
That is, for every fixed $\theta$, the projection $R_Hf_{\theta}$ defined by $R_Hf_{\theta}(s)=Rf(s,\theta)$ (which is a one dimensional cumulative density function) can be approximated by the empirical distribution function of the projection of $X$ to the line $\{s\theta : s\in\mathbb{R}\}$ \cite{rad2}. Similarly for the spherical Radon transform we have
\begin{equation}
\begin{split}
\label{equ2.3}
R_Mf(s,x)=\int_{0}^{s}r^{n-1}\int_{S^{n-1}}f(x+r\theta)\mathrm{d}\theta\mathrm{d}r\approx\sum_{i\in\{(x_i-x)^2<s^2\}}1.
\end{split}
\end{equation}
So by counting the number of observations which lie in a set of half spaces or spheres we can approximate Radon transforms of the density $f$. See figure \ref{fig1}. 

We aim to discretize the Radon transform operators above to a pixel grid (assuming a piecewise constant density as is commonly applied in image reconstruction and inverse problems) and solve the resulting sparse system of linear equations to recover the density. Here we will make use of widely applied regularization techniques (such as TV) to combat the noise due to finite sampling, and apply known automated methods in tomography to choose the regularization parameter (such as the Generalized Cross Validation (GCV)). First we provide theoretical performance garuantees for our method, in the next section.


\section{Error estimates}
\label{error1}
Here we give error estimates for our method for density estimation combining results on the stability of inverse Radon transforms and the expected error of empirical distribution functions. The next result is a statement of the Dvoretzky-Kiefer-Wolfowitz (DKW) inequality \cite{DKW} which explains the error in an approximation of a univariate density from its empirical distribution function.

We now have our main theorem which gives bounds for the error in a half space Radon transform approximation from a full set of projected one dimensional empirical distribution functions.
\begin{theorem}
\label{errthm}
Let $x_1,\ldots,x_m$ be continuous, independant and identically distributed random variables in $\Omega^n$ with probability density function $f \in C^{\infty}_0(\Omega^n)$ and let $$g_{m}(s,\theta)=\frac{1}{m}\sum_{i\in\{x_i\cdot\theta\leq s\}}1$$ be an approximation to the half space Radon transform $g=R_Hf$. Let $K$ projection directions $\{\theta_j\}_{j=1}^K$ be uniformly spread over $S^{n-1}$. Then
\begin{equation}
\|g-g_{m}\|^2_{L^2(Z)}\leq w_{1,n-1}\frac{\log\frac{2K}{p}}{m}+\epsilon(K)
\end{equation}
with probability $1-p$ for any $0\leq p\leq 1$, where $g_{\theta}(s)=g(s,\theta)$ and $g_{m,\theta}(s)=g^{\epsilon}(s,\theta)$ and 
$$\epsilon(K)=\left|\|g-g_m\|^2_{L^2(Z)}-\sum_{j=1}^Kw_j\|g_{\theta_j}-g_{m,\theta_j}\|^2_{L^2(\mathbb{R})}\right|$$
is the error term in approximating the integral over $S^{n-1}$ by a Reimann sum, where the area elements $w_j$ are such that  $\sum_{j=1}^Kw_j=w_{1,n-1}$ and $w_{1,n-1}$ is the surface area of $S^{n-1}$.
\begin{proof}
The proof is deferred to section \ref{sec:proof}.
\end{proof}
\end{theorem}
\begin{corollary}
\label{errcorr}
Let $x_1,\ldots,x_m$ be continuous, independant and identically distributed random variables in $\Omega^n$ with probability density function $f \in C^{\infty}_0(\Omega^n)$ and let $$g_m(s,\theta)=\frac{1}{m}\sum_{i\in\{x_i\cdot\theta\leq s\}}1$$ be an approximation to the half space Radon transform $g=R_Hf$. Then the worst case error in a reconstruction $f_m$ of $f$ from $g_m$ is
\begin{equation}
\|f-f_m\|_{L^2(\mathbb{R}^n)}\leq c(n)\left(w_{1,n-1}\frac{\log\frac{2K}{p}}{m}+\epsilon(K)\right)^{\frac{1}{2(n+2)}}\rho^{1-1/(n+2)},
\end{equation}
which holds with probability $1-p$. Here $K$ is the number of projections (which are assumed to be uniformly distributed on the sphere), $\|f\|_{H^{\frac{1}{2}}}\leq\rho$, and $\epsilon(K)$ is an error term in approximating the integral over $S^{n-1}$ by a finite sum.
\begin{proof}
The result follows immediately from Theorem \ref{errthm} and Theorem \ref{Rerrthm1}.
\end{proof}
\end{corollary} 

Regarding the above results and the error term $\epsilon(K)$, the idea here is that we can pick $K$ large enough (indeed we have no restriction (in theory) on the number of projections taken) to imply a stable inversion (or such that $\epsilon(K)\approx 0$). That is, we can choose $K$ so that the inverse problem of reconstructing $f$ from $R_Hf$ (known for a finite number of projections) is (effectively) mildly ill--posed, and such that there are no image artefacts, in the sense that the singularities of $f$ are stably recovered in all directions \cite{quinto2,quinto3} (hence the choice of uniformity on $S^{n-1}$ for the set of directions $\{\theta_j\}_{j=1}^K$). 

Let $K=k^{n-1}$ (we expect the number of projections needed to adequately cover $S^{n-1}$ to increase exponentially with $n$) be large enough so that the error $\epsilon(K)\approx 0$ is negligible. For example we could fix $k=360$ for an angular step size of $1^{\circ}$ on $[0,2\pi]\cong S^1$. Then
\begin{equation}
\begin{split}
\|g-g_m\|^2_{L^2(Z)}\leq w_{1,n-1}\frac{\log\frac{2k^{n-1}}{p}}{m}=\frac{w_{1,n-1}}{m}\left((n-1)\log(k)+\log\left(\frac{2}{p}\right)\right)
\end{split}
\end{equation}
with probability $1-p$. The bound above explains how we would expect the error in our Radon transform approximation to increase with $n$, $m$ and the confidence $1-p$. 
For fixed $n$, $k$ and $p$ we expect the convergence rate to the solution of $(1/m)^{\frac{1}{2(n+2)}}$ with $m$. 


The next theorem describes the least squares error expected in a spherical Radon transform approximation using empirical distribution functions.
\begin{theorem}
\label{sphthm}
Let $x_1,\ldots,x_m$ be continuous, independant and identically distributed random variables in $\Omega^n$ with probability density function $f \in C^{\infty}_0(\Omega^n)$ and let $$g_m(s,x)=\frac{1}{m}\sum_{i\in\{(x_i-x)^2\leq s^2\}}1$$ be an approximation to the spherical Radon transform $g=R_Mf$. Then
\begin{equation}
\frac{1}{K}\sum_{j=1}^{K}\|g_{x_j}-g_{m,x_j}\|^2_{L^2(\mathbb{R}^+)}\leq\frac{1}{m}\left(\log(K)+\log\left(\frac{2}{p}\right)\right)
\end{equation}
with probability $1-p$, for any finite set $\{x_j\}_{j=1}^K$ of circle centers. Here $g_{x_j}(s)=g(s,x_j)$ and $g_{m,x_j}(s)=g_m(s,x_j)$.
\begin{proof}
The proof is deferred to section \ref{sec:proof}. See Theorem \ref{sphthm1}.
\end{proof}
\end{theorem}
So we would expect the same rate of decrease in the error in a spherical Radon transform approximation as with half space data. The spherical Radon transform is mildly ill--posed and to the same degree as the half space transform if enough spherical centers $x$ and radii $s$ are known (this can be shown using the theory presented in \cite{quinto}). So we would expect the same amplification in the error in our solution. In general, if we approximate Radon transform projections (e.g. ellipsoidal, hyperbolic) as empirical distribution functions, we can expect the same rate of convergence, provided that the Radon transforms used have a stable inverse. 


\section{Method and results}
\label{method}
Here we present practical methods for density estimation in low dimensions by a discrete inversion of Radon transforms. We test our approach on a number of synthetic densities in two dimensions and give a side by side comparison with a kernel estimator. Later in section \ref{manifolds} we expand our approach to density estimation on manifolds and provide an error analysis of the reconstructions of densities on manifold patches of varying levels of curvature and radius.

There are two algorithms presented in this section, which we will now introduce. Algorithm \ref{alg1} implements density estimation in low dimensions by an explicit, discrete inversion of Radon transforms with regularization. Algorithm \ref{alg3} applies Algorithm \ref{alg1} to local tangent spaces of manifolds to reconstruct densities on low dimensional manifolds embedded in $\mathbb{R}^d$. 
\subsection{Density estimation in low dimensions}
To implement our density estimation method in low dimensional ($n\leq 3$) Euclidean space we will approximate Radon transform operators as discrete sums over a uniform grid (which would typically represent image pixels in image reconstruction). Specifically we aim to minimize the following functional
\begin{equation}
\label{lls}
\|Rv-b\|_2^2+\lambda^2\mathcal{G}(v),
\end{equation}
where $R$ is the discrete form of the Radon transform (either spherical or half space) and $\mathcal{G}$ is a regularization penalty with regularization parameter $\lambda$. For example $\mathcal{G}(v)=\|v\|_2^2$ (Tikhonov) or $\mathcal{G}(v)=\sum_i|v_{i}-v_{i-1}|$ (total variation). To solve the least squares problem (\ref{lls}) we apply the iterative solvers of \cite{irtools}, and choose the regularization parameter $\lambda$ via a Generalized Cross Validation (GCV) method \cite[page 95]{hansen1}. To estimate densities in low dimensions we implement Algorithm \ref{alg1}.
\begin{algorithm}[!h]
\KwResult{Density estimation by inverse Radon transform}
\hspace*{\algorithmicindent} \textbf{Input:} A point cloud $X\subset\mathbb{R}^n$, a pixel grid which covers $X$, with pixel centers $\{p_i\}_{i=1}^{l_1}\subset\mathbb{R}^n$, a set of balls $\{s_j\}_{j=1}^{l_2}$ and half spaces $\{h_j\}_{j=1}^{l_2}$ in $\mathbb{R}^n$\\
 \textbf{Output:} A vector of density values $v$, where $v_i=f(p_i)\in \mathbb{R}$ for $1\leq i\leq l_1$\\
\textbf{Procedure:}\\
Set the transform ($T$), either spherical ($T=T_1$) or half space ($T=T_2$)\\
Initialize a zero sparse matrix $R$ with $l_1$ columns and $l_2$ rows, and a zero vector $b$ with $l_2$ entries\\
\eIf{$T=T_1$}{ \For{$j\leftarrow 1$ \KwTo $l_2$}{
Find the $p_i$ such that $p_i\in s_j$ and set $R_{ji}=1\iff p_i\in s_j$\\
Count the number of observations in $s_j$ and set $b_j=\left|X\cap s_j\right|$
}
}{ \For{$j\leftarrow 1$ \KwTo $l_2$}{
Find the $p_i$ such that $p_i\in h_j$ and set $R_{ji}=1\iff p_i\in h_j$\\
Count the number of observations in $h_j$ and set $b_j=\left|X\cap h_j\right|$
}
}
Set the regularization function $\mathcal{G}$ and find
$$\argmin_v\|Rv-b\|_2^2+\lambda^2\mathcal{G}(v)$$
using an iterative solver, choosing $\lambda$ by a GCV technique\\
Output $v$\;
 \caption{Density estimation by inverse Radon transform}
 \label{alg1}
\end{algorithm}

\begin{figure}[!h]
\centering
\begin{subfigure}{0.35\textwidth}
\includegraphics[width=1\linewidth, height=4.8cm]{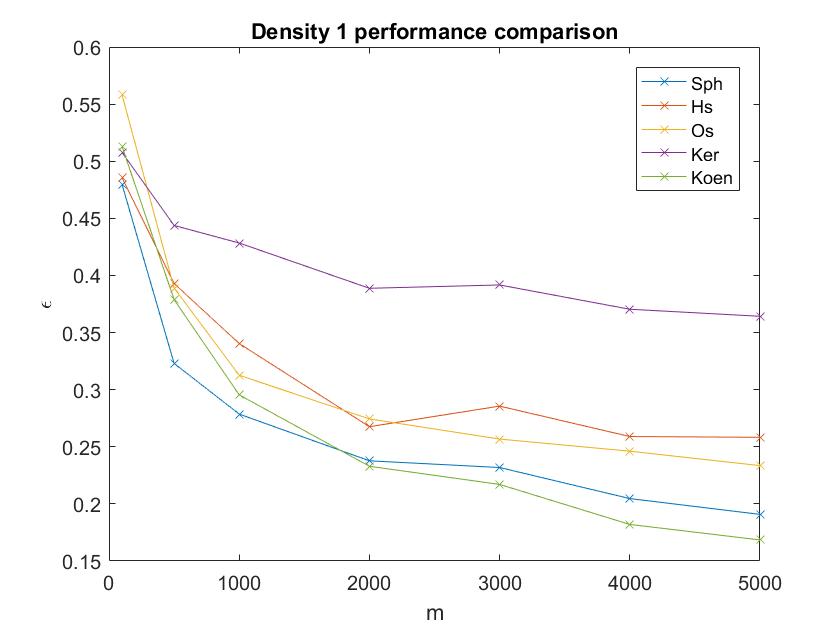} 
\end{subfigure}
\begin{subfigure}{0.35\textwidth}
\includegraphics[width=1\linewidth, height=4.8cm]{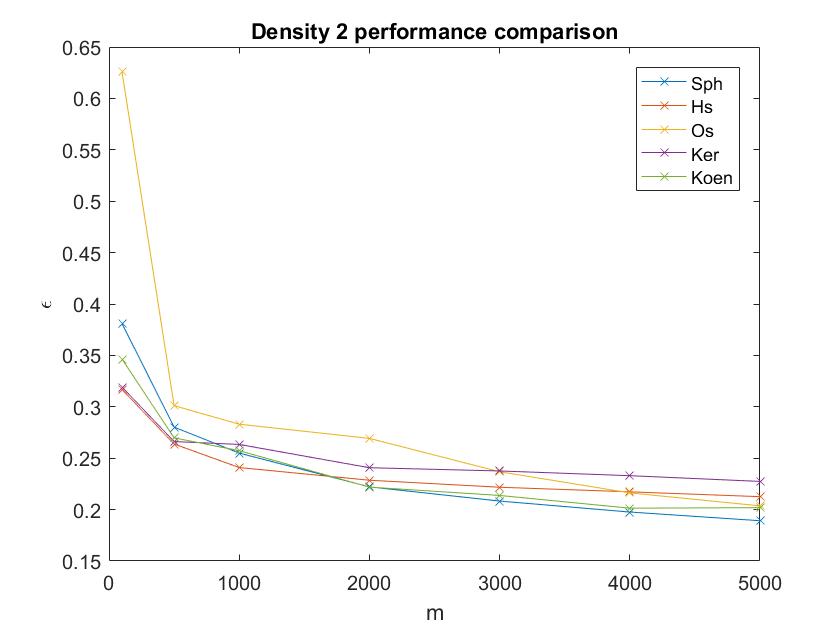}
\end{subfigure}
\begin{subfigure}{0.35\textwidth}
\includegraphics[width=1\linewidth, height=4.8cm]{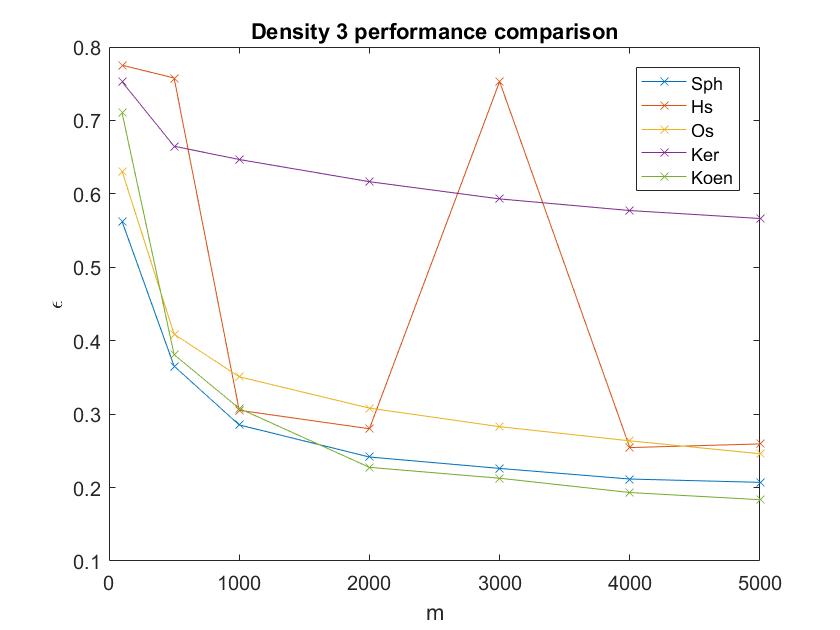} 
\end{subfigure}
\begin{subfigure}{0.35\textwidth}
\includegraphics[width=1\linewidth, height=4.8cm]{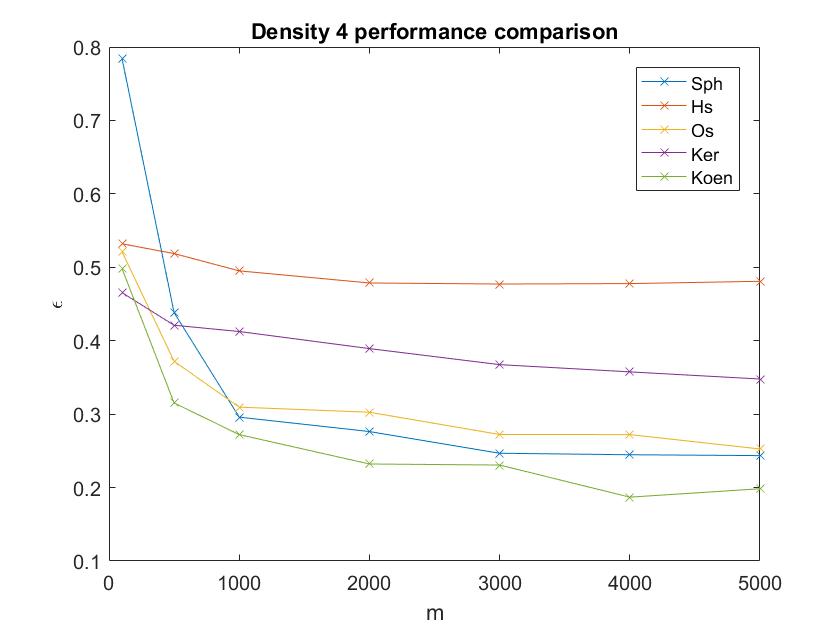}
\end{subfigure}
\caption{Comparisons of errors $\epsilon$ in density reconstructions 1--4, using the methods Sph, Hs, Os, Ker and Koen.}
\label{figC}
\end{figure}

To show the effectiveness of our method we give an error and image reconstruction comparison with a collection of density estimation methods from the literature using TV and Radon transforms. The first is that of Koenker {\it{et. al.}} \cite{TV1D1} (denoted by Koen), which finds the density maximizing the log--likelyhood function with TV penalties. In \cite{TV1D1} the selection of $\lambda$ is not automated and they set $\lambda=2$. We could employ the use (ordinary) cross validation (CV) methods here to choose $\lambda$. However this would require an impractically long run time, as is discussed in \cite{mohler}. Further, the use of GCV methods would be ill--advised due to the arguments given in \cite[page 11]{TV1D}. Hence, for every density and sample size $m$ considered, we report the best result (in terms of $\epsilon$) over three levels of smoothing, namely $\lambda=0.5,1,2$. The second point of comparison is a Gaussian kernel estimate with a uniform badwidth kernel (denoted by Ker), with the bandwidth chosen to be optimal for Gaussian densities \cite{gdense}. 
The third point of comparison considered in the main text is that of O'Sullivan {\it{et. al.}} \cite{rad2} (denoted by Os). Here we use the truncated Fourier series approximations to the Radon projections $g_{\theta}$ given on page 518 of \cite{rad2}. Such approximation ideas are shown to the satisfy the performance guarantees provided in \cite{rad2}. We set
\begin{equation}
\label{Osapp}
\hat{g}_\theta(\xi)=\text{rect}\left(\frac{\xi}{2h_m}\right)\cdot\frac{1}{m}\sum_{i=1}^m\exp(i\xi\theta\cdot x_i),
\end{equation}
where $\hat{g}_{\theta}$ denotes the Fourier transform of $g_{\theta}$, $\text{rect}$ is the rectangle function and $h_m$ is a chosen bandwidth parameter. This approximation is equivalent to a 1-D kernel estimate using sinc functions. That is $g_{\theta}(x)=\frac{1}{m}\sum_{i=1}^m\sinc(\frac{x-\theta\cdot x_i}{h_{m}})$. In \cite[page 518]{rad2} a choice for $h_m$ is suggested based on the size of the Sobolev norms $\|g_{\theta}\|_{H^k}$. However, as $\|g_{\theta}\|_{H^k}$ is not known a-priori, we report the best result (in terms of $\epsilon$) over three levels of smoothing $h_m=0.5,1,2$ (scaling the sinc functions to half and twice the pixel length of 1) for every density and sample size $m$ considered. We calculate (\ref{Osapp}) for $\theta = \{\frac{\pi i}{180} : 0\leq i \leq 179\}$ and then reconstruct the density using the filtered back-projection formula, as in \cite{rad2}. In the appendix, we also present the density estimation errors using TV and Poisson image denoising approaches, and using the Split Bregman techniques of \cite{mohler}, which offer a faster implementation of the density estimation method suggested by Koenker \cite{TV1D1}. The results using these methods are left to the appendix (see figure \ref{figCS}) as they were found to produce a significantly higher error than those considered in figure \ref{figC}. Hence we focus on the more competitive approaches in the main text. It seems, based on the current experiments, that a direct density estimation using TV (by minimizing the least squares error as in \cite{ROF}) is ill advised for the types of densities considered, and we see a better performance with a density estimation in the Radon domain combined with a TV regularizer in the inversion. 

Let $v$ be the vectorized form of the original density, and let $v_m$ denote an approximation. Then we measure the error by
$
\label{errmeas}
\epsilon=\frac{\|v-v_m\|_2}{\|v\|_2}.
$
We consider four test densities, whose probabilty density functions are displayed in the left hand of figures \ref{fig4}--\ref{fig8}. See figure \ref{fig3.3} for a larger view and description of the ground truth densities. We give the estimation errors for $m\in\{100,500,1000,2000,3000,4000,5000\}$ and present image reconstruction comparisons for the two most competitive methods (with the least error) for $m=2000$. All densities considered are piecewise smooth, some with more edges and more smooth features than others. We use these as example densities to test the effectiveness of our reconstruction method with TV for recovering the jump singularities in the density while also preserving smoothness. Each of the test densities displayed are to be reconstructed on a $100$--$100$ pixel grid (the density is assumed to be piecewise constant with $100^2$ densities, the matrix $R$ has $100^2$ columns). For half space Radon transform data, we sampled $s$ and $\theta$ as $s\in \{-50+i : 0\leq i\leq 100 \}$ and $\theta = \{\frac{\pi i}{180} : 0\leq i \leq 179\}$ to be perfectly uniform on $Z$, as is most optimal by Theorem \ref{sample}, and in line with the theory presented in the previous section. For spherical Radon transform data, we sample circle centers $x\in \{(0.5+i,0.5+j) : 0\leq i,j\leq 99\}$ and radii $s\in \{4+i : 0\leq i\leq 16\}$. Counting observations in circles of smaller radii ($s=4$), while being a more noisy approximation, helps to highlight finer features in the density, whereas the larger radii ($r=20$) offer less precise information but with less noise. We find that using a range of radii is most optimal and overall gives better performance than using half space data (see appendix \ref{appC} for a full comparison), which is what we would expect as spherical Radon transform data is overdetermined.

We find that a density reconstruction using Sph and Koen gave the best and most consistent performance in terms of $\epsilon$ for all densities considered, and for varying sample sizes $m$. See figure \ref{figC} for a comparison of the errors using our methods (spherical and half space, denoted as Sph and Hs respectively), and the methods Os \cite{rad2}, Koen \cite{TV1D1} and Ker.  For a full comparison of spherical vs half space consult appendix \ref{appC}. To generate the Gaussian mixture in density 1 we randomly (uniformly) selected 100 Gaussian centers (means) on a 1--100 meshgrid, with constant covariance, and took the sum. To see how the methods considered compare over a variety of Gaussian mixture densities, see table \ref{T1}, where we have given the mean and standard deviations of the errors in a reconstruction of density 1 with 20 random Gaussian mean sets and a constant sample size $m=1000$. As we see an occasional failure in the reconstruction of density 3 using Hs (we see a spike in the error in the red curve in the bottom left hand of figure \ref{figC} at $m=3000$), we also give the mean and standard deviation reconstruction errors from 20 random draws of density 3, with $m=1000$ samples. See table \ref{T2}. We see further evidence of a more robust and accurate performance using Sph and Koen when compared to the other methods. Although as suspected, the range and mean of errors is far greater using Hs, so there are cases when Algorithm \ref{alg1} may break down using Hs techniques.
\begin{table}[H]
\centering
\begin{tabular}{| c | c | c | c | c | c | c | c |}
\hline
$\epsilon$ & Sph  & Hs   & Os   & Ker  & Koen \\ \hline
$\mu$      & 0.29 & 0.32 & 0.32 & 0.40  & 0.29 \\ 
$\sigma$   & 0.03 & 0.03 & 0.02 & 0.04 & 0.02 \\ \hline
\end{tabular}
\caption{Mean and standard deviation errors from 20 random draws ($m=1000$) from 20 different Gaussian mixtures with randomly generated means (similar to density 1).}
\label{T1}
\end{table}
\begin{table}[H]
\centering
\begin{tabular}{| c | c | c | c | c | c | c | c |}
\hline
$\epsilon$ & Sph  & Hs   & Os   & Ker  & Koen \\ \hline
$\mu$      & 0.32 & 0.48 & 0.34 & 0.65 & 0.28 \\ 
$\sigma$   & 0.01 & 0.20  & 0.01 & 0.01 & 0.01 \\ \hline
\end{tabular}
\caption{Mean and standard deviation errors from 20 random draws ($m=1000$) of density 3.}
\label{T2}
\end{table}

Refer to figures \ref{fig4}--\ref{fig8}, where we have presented density reconstruction comparisons (as 2-D images) for the two most competitive methods (with the least error) on each density 1--4, for $m=2000$. In each case the most competitive method to Sph, was Koen. We see a similar image reconstruction quality and error level when we compare Sph and Koen across densities and sample size. Although we see a greater image quality using Koen for densities 3 and 4, and conversely using Sph for densities 1 and 2. It is noted that the choice of smoothing parameter is not automated using Koen, and in \cite{TV1D1} no such method is suggested. In \cite{mohler} Mohler {\it{et. al.}} provide a faster implementation of the Koen method using Split Bregman techniques, which allows for an automated choice of $\lambda$ by 10-fold CV. The increased speed of implementation seems to be at a detrement to the error however, as can been from the results presented in \ref{addtab}, and the errors when using the faster implementation of Koen are not competitive with Sph. As the choice of the smoothing $\lambda$ is one of the most crucial aspects in density estimation (which we can fully automate using GCV with spherical Radon transforms, given the linear inversion), and there is no known automated method for choosing $\lambda$ using the more accurate but less efficient method of \cite{TV1D1}, on this basis, the above discussions and the experiments presented here, we conclude that a spherical Radon approach is optimal for the types of densities considered. 

\begin{figure}[!h]
\begin{subfigure}{0.32\textwidth}
\includegraphics[width=0.9\linewidth, height=4cm]{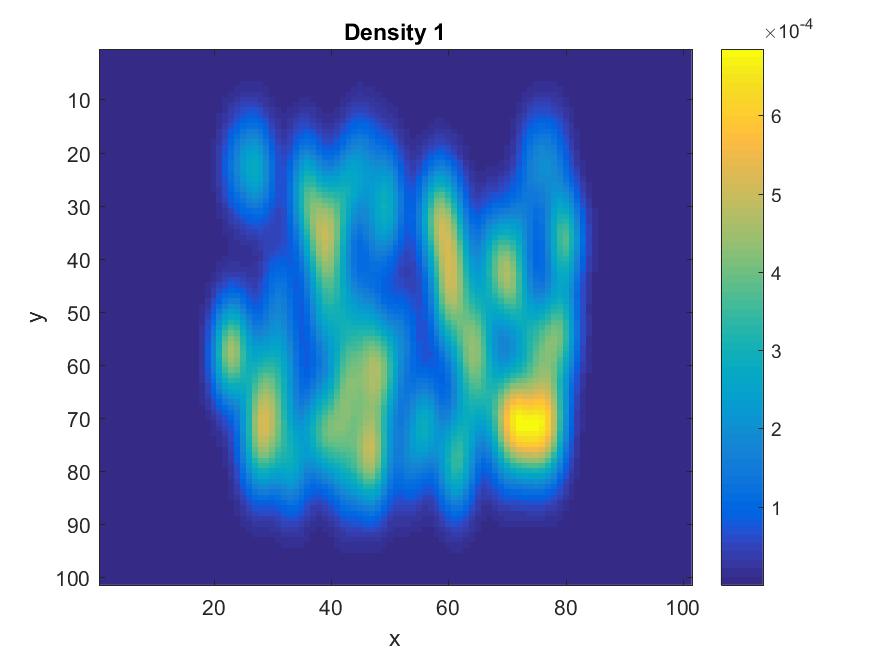}
\end{subfigure}
\begin{subfigure}{0.32\textwidth}
\includegraphics[width=0.9\linewidth, height=4cm]{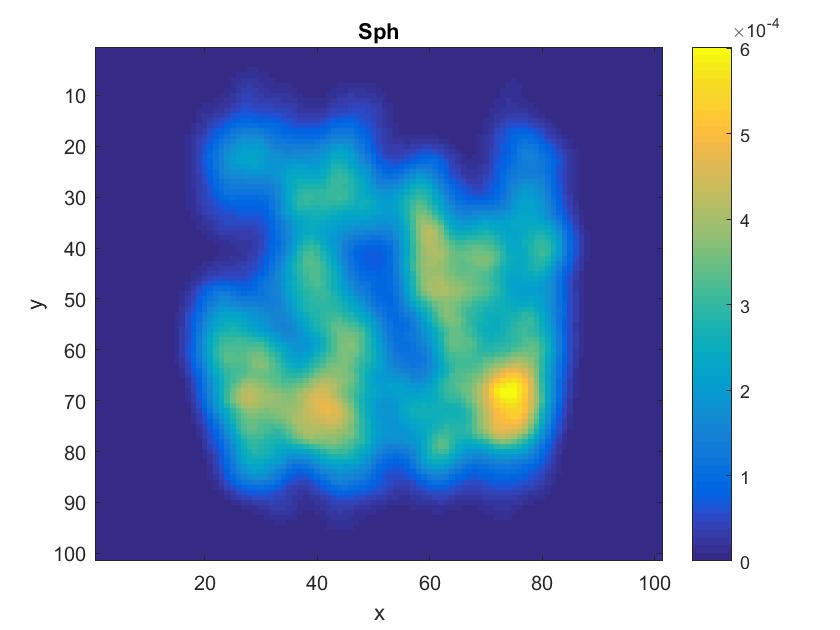} 
\end{subfigure}
\begin{subfigure}{0.32\textwidth}
\includegraphics[width=0.9\linewidth, height=4cm]{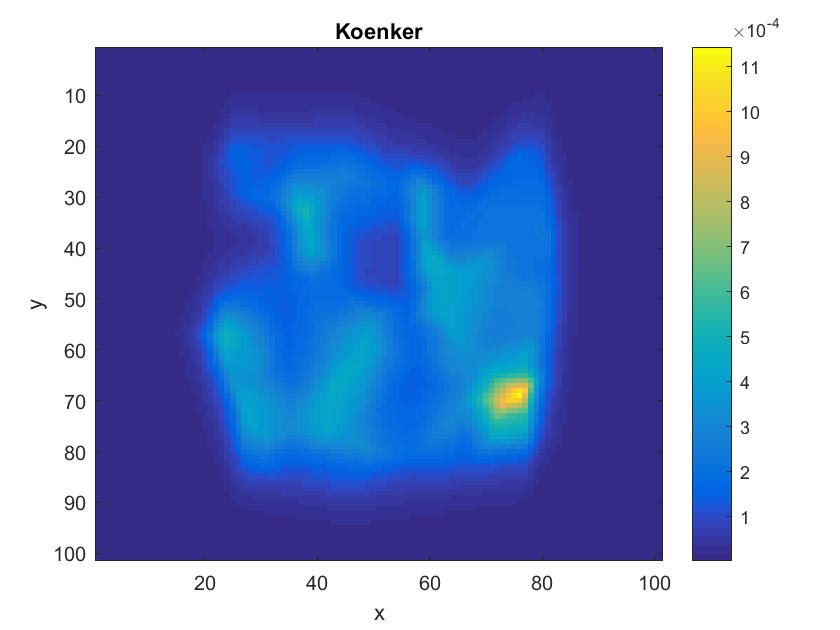}
\end{subfigure}
\caption{From left to right, Density 1, a reconstruction using Sph, and Koen.}
\end{figure}
\begin{figure}[!h]
\begin{subfigure}{0.32\textwidth}
\includegraphics[width=0.9\linewidth, height=4cm]{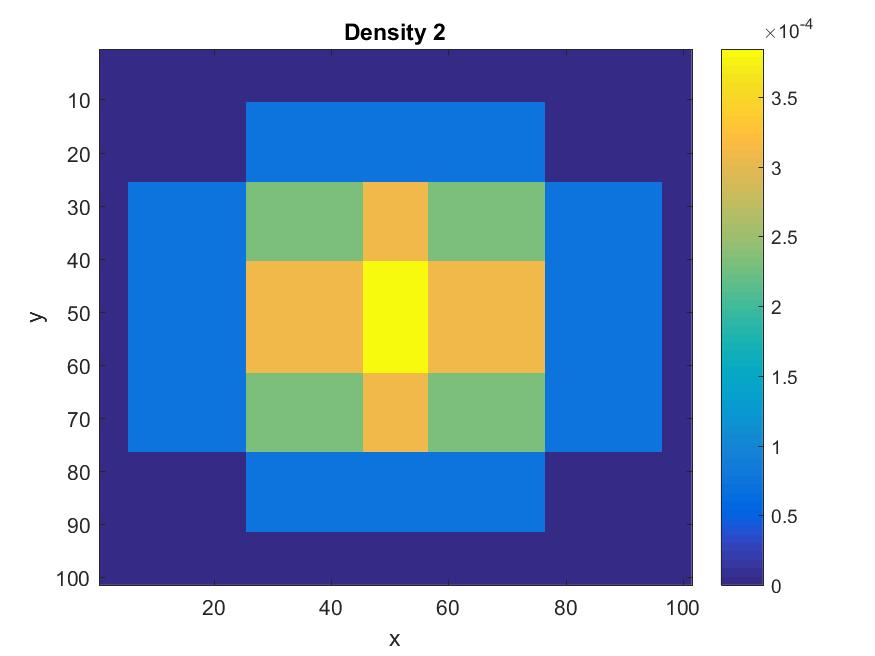}
\end{subfigure}
\begin{subfigure}{0.32\textwidth}
\includegraphics[width=0.9\linewidth, height=4cm]{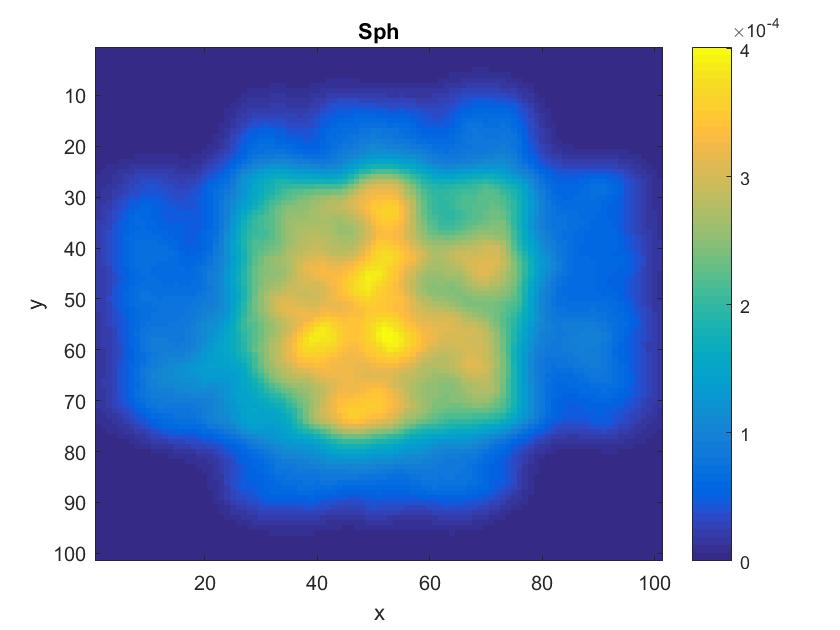} 
\end{subfigure}
\begin{subfigure}{0.32\textwidth}
\includegraphics[width=0.9\linewidth, height=4cm]{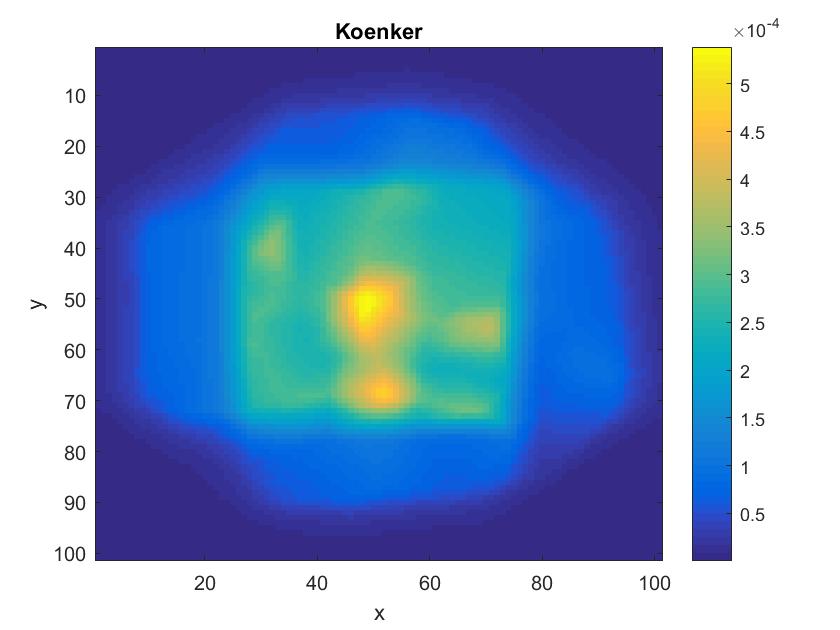}
\end{subfigure}
\caption{From left to right, Density 2, a reconstruction using Sph, and Koen.}
\label{fig8}
\end{figure}
\begin{figure}[!h]
\begin{subfigure}{0.32\textwidth}
\includegraphics[width=0.9\linewidth, height=4cm]{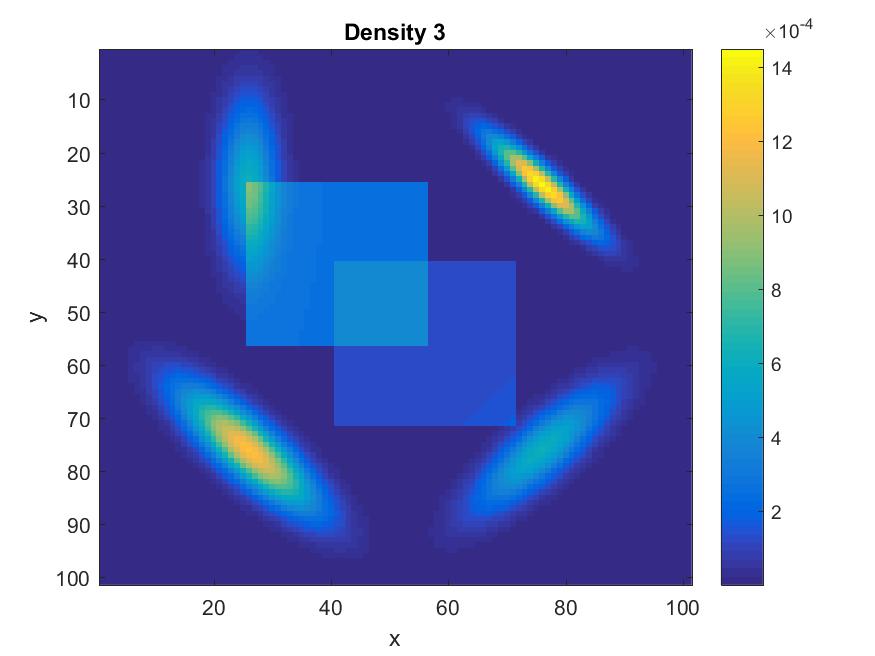}
\end{subfigure}
\begin{subfigure}{0.32\textwidth}
\includegraphics[width=0.9\linewidth, height=4cm]{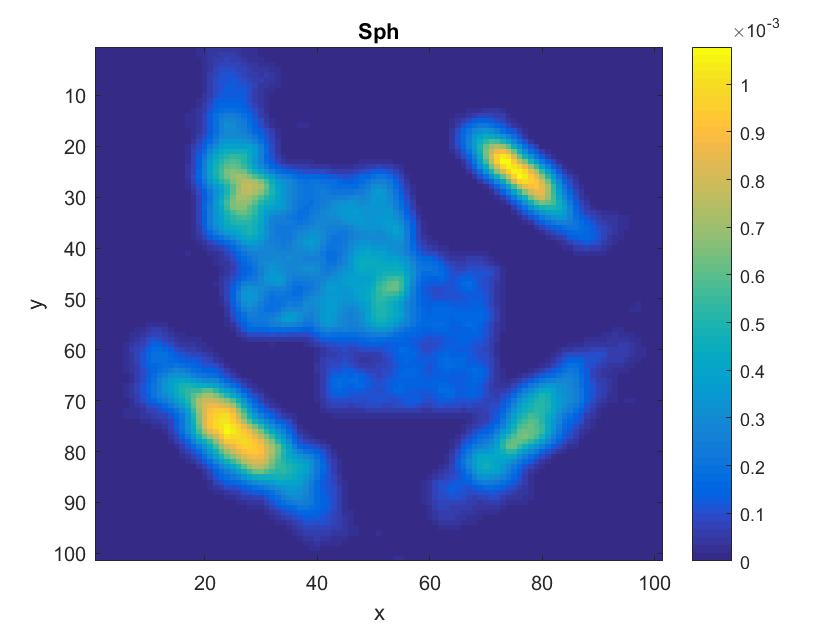} 
\end{subfigure}
\begin{subfigure}{0.32\textwidth}
\includegraphics[width=0.9\linewidth, height=4cm]{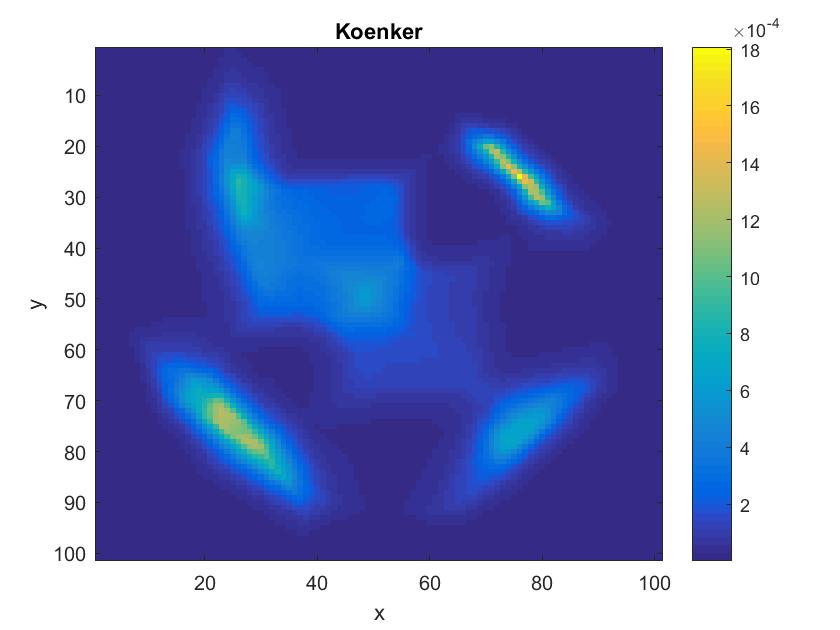}
\end{subfigure}
\caption{From left to right, Density 3, a reconstruction using Sph, and Koen.}
\label{fig4}
\end{figure}
\begin{figure}[!h]
\begin{subfigure}{0.32\textwidth}
\includegraphics[width=0.9\linewidth, height=4cm]{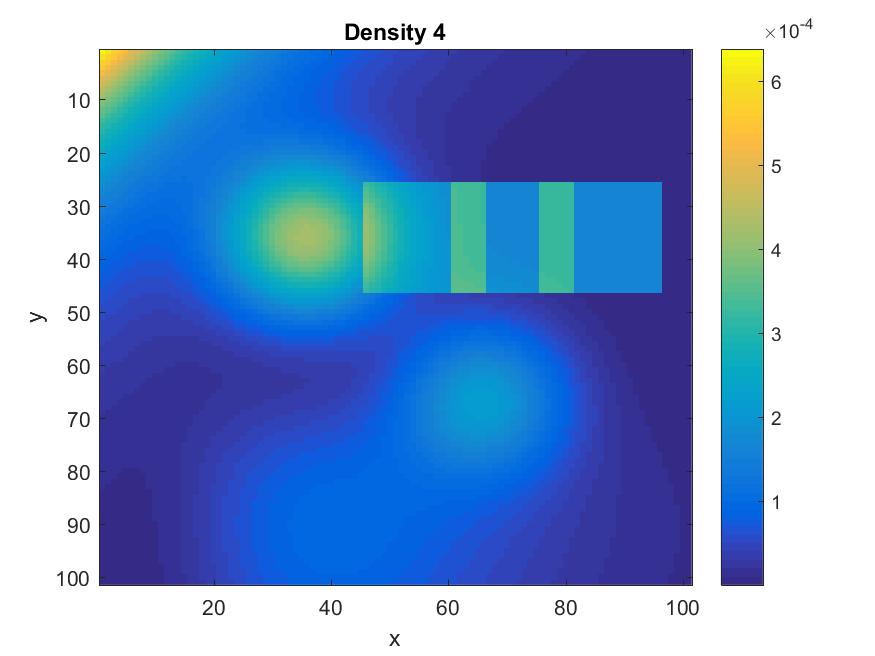}
\end{subfigure}
\begin{subfigure}{0.32\textwidth}
\includegraphics[width=0.9\linewidth, height=4cm]{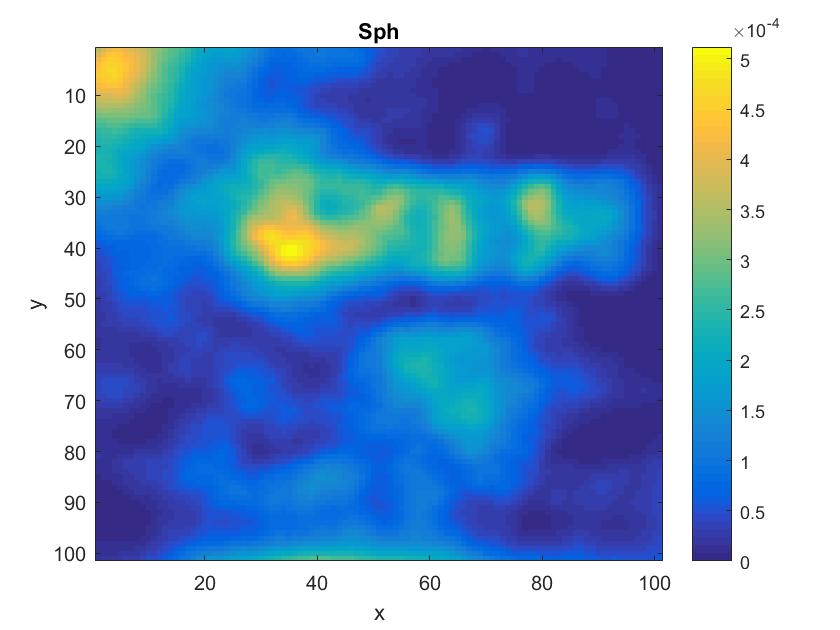} 
\end{subfigure}
\begin{subfigure}{0.32\textwidth}
\includegraphics[width=0.9\linewidth, height=4cm]{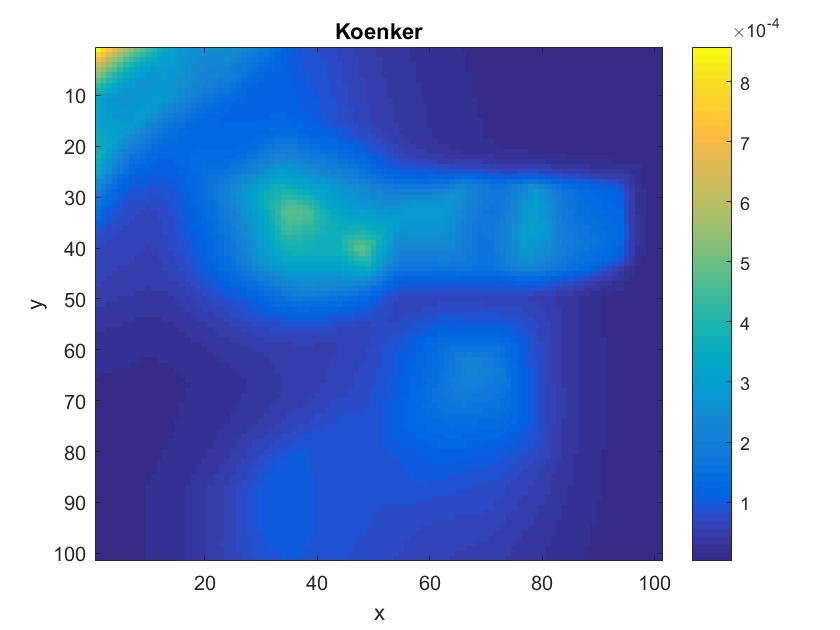}
\end{subfigure}
\caption{From left to right, Density 4, a reconstruction using Sph, and Koen.}
\label{fig6}
\end{figure}
\subsection{Density estimation on low dimensional manifolds}
\label{manifolds}
Here we extend our previous results to density estimation on low dimensional manifolds embedded in $\mathbb{R}^d$, where we will incorporate our previous reconstruction methods and theory, the theorems of \cite{fefferman} and ideas in manifold learning. First some definitions. Let $X$ be a metric space. Throughout this section we denote $B^X_r(x)$ to be the ball of radius $r$ centerd at $x$ and $B^n_r(x)=B^X_r(x)$ for $X=\mathbb{R}^n$.

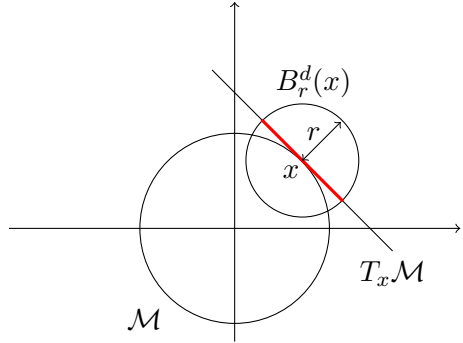
\begin{figure}[!h]
\centering
\begin{tikzpicture}[scale=3]
\draw [->] (0,-0.5)--(0,1);
\draw [->] (-1,0)--(1,0);
\draw (-0.1,0.7)--(0.7,-0.1);
\draw (0,0) circle (0.42);
\node at (0.25,0.25) {$x$};
\node at (0.7,-0.2) {$T_x\mathcal{M}$};
\node at (-0.4,-0.4) {$\mathcal{M}$};
\draw (0.3,0.3) circle (0.25);
\node at (0.35,0.65) {$B^d_r(x)$};
\draw [red,very thick] (0.12,0.48)--(0.48,0.12);
\draw [<->] (0.3,0.3)--(0.47,0.47);
\node at (0.35,0.42) {$r$};
\end{tikzpicture}
\caption{A manifold $\mathcal{M}=S^1$ in $\mathbb{R}^d$, where $d=2$, and a visualisation of $B^{T_x\mathcal{M}}_r(x)$, which in this case is a line segment (highlighted in red).}
\label{fig2.2}
\end{figure}

We now have the theorem from \cite{fefferman} which gives bounds for the error in a local tangent space approximation of a manifold $\mathcal{M}$ in terms of its principle curvatures.
\begin{theorem}
\label{Tspace}
Let $\mathcal{M}$ be a submanifold of $\mathbb{R}^d$ whose principle curvatures are bounded by $\kappa>0$. Then
\begin{equation}
d_H(B^{\mathcal{M}}_r(x),B^{T_x\mathcal{M}}_r(x))\leq \frac{1}{2}\kappa r^2
\end{equation}
for any $x\in \mathcal{M}$ and $r\geq 0$, where $B^{T_x\mathcal{M}}_r(x)=B^d_r(x)\cap T_x\mathcal{M}$ and $d_H$ denotes the Hausdorff distance.
\end{theorem}
See figure \ref{fig2.2} for a visualisation of a set $B^{T_x\mathcal{M}}_r(x)$.

So far we have dealt with density estimation in low dimensions, where the true data dimension is known ($n=2$ in the above simulations), in the sense that the observation sets $X$ lie on 2--manifolds in $\mathbb{R}^2$. In many machine learning applications however, the true data dimension $n$ is hidden and and we are given a set of observations $x'_i\in\mathbb{R}^d$, where $n<d$. It is often assumed that the observations lie on a manifold $\mathcal{M}$ which is embedded in $\mathbb{R}^d$. More precisely for every $1\leq i\leq m$, $x'_i=\phi(x_i)$, where $\phi : \mathcal{M}\to\mathbb{R}^d$ is an embedding. This is the idea behind the field of ``manifold learning".

Let $\mathcal{M}\subset\mathbb{R}^n$ be an $n$--manifold and let our observation set $X=\{x_1,\ldots,x_m\}\subset \mathcal{M}$. Let us assume that the data we are given is $X'=\{\phi(x_1),\ldots,\phi(x_m)\}$, where $\phi : \mathcal{M}\to\mathbb{R}^d$ is an embedding in the general topological sense (a homeomorphism onto its image), with $d>n$. Let $\mathcal{M}$ be described by a collection of charts $\{U_l,\varphi_l\}$, where $U_l$ are open sets which cover $\mathcal{M}$ and each $\varphi_l :U_l\to\Omega^n$ is a homeomorphism of $U_l$ to the unit ball in $\mathbb{R}^n$. Let $f :\mathcal{M}\to\mathbb{R}$ be a density on $\mathcal{M}$. Then $f'=f \circ \phi^{-1} : \phi(\mathcal{M})\to\mathbb{R}$ is the representation of $f$ in $d$ dimensions. If the charting $\{U_l,\varphi_l\}$ is known, then we can recover densities on $\mathcal{M}$ using Radon transforms. 

More precisely, for every $x'\in \phi(\mathcal{M})$, there exists an $l$ such that $x'\in U_l$, and the function $h=f'\circ \varphi^{-1}_l :\Omega^n\to \mathbb{R}$ can be reconstructed explicitly from its spherical or half space Radon transform. Once the function $h$ is known we can compose it with $\varphi_l$ to recover the corresponding density value for $x'$. For a finite sample $X'$ on $\mathcal{M}$ the error in the reconstruction of a patch $U_l$ is described by Corollary \ref{errcorr}. 
This is under the assumption of course that the charting of $\mathcal{M}$ is known, which is not often the case. We aim to tackle this problem in this section.
\begin{table}[H]
\centering
	\begin{tabular}{| c | c | c | c | c | c | c | c |}
	\hline
		$\epsilon$ & $\kappa=0.01$ & $\kappa=0.05$ & $\kappa=0.1$     \\ \hline
$r=5$ &       $0.36$ &	$0.39$ &	$0.38$   \\ 
$r=10$ &       $0.22$ &	$0.27$ &	$0.36$   \\ 
$r=20$ &       $0.17$ &	$0.23$ &	$0.34$    \\ 
\hline
\end{tabular}
\caption{Relative errors in patch density reconstructions for varying levels of curvature ($\kappa$) and patch radius ($r$). The reconstructions were produced using Algorithm \ref{alg3} with the method set to spherical and the variance percentage set to $p=90\%$.}
\label{tabEmb}
\end{table}
\subsubsection{Local tangent space approximations}
Let $\{U_l,\varphi_l\}$ be an atlas of an $n$--manifold $\mathcal{M}$ embedded in $\mathbb{R}^d$, with $n<d$, and let each $U_l=B_r^{\mathcal{M}}(x)$ be a neighborhood of a point $x\in\mathcal{M}$ for some $r$. By Theorem \ref{Tspace}, we can approximate $U_l$ by the tangent patch $B^{T_x\mathcal{M}}_r(x)=T_x\mathcal{M}\cap B^d_r(x)$, provided the principle curvatures of $\mathcal{M}$ are low or $r$ is small enough. This is an idea used in techniques in manifold learning and non-linear dimensionality reduction \cite{fefferman,LTSA}. For finite data $X'\subset\mathcal{M}$ we will take subsamples of the data $S'=X'\cap B^d_r(x)$ (to represent a sample on $U_l$), for $r$ small enough, and project onto the first $n$ principle components. This approximates the projection of $S'$ onto $T_x\mathcal{M}$ (or the application of the chart $\varphi_l :U_l\to\Omega^n$). So here the atlas $\{U_l,\varphi_l\}$ which describes $\mathcal{M}$ is such that each $U_l$ is approximately flat. After we have transformed the samples $S'$ to principle component space, and let us denote the projected subsample by $S$, we can reconstruct the density from which $S$ is drawn (this would be $h$, using the notation above) by applying Algorithm \ref{alg1}. The reconstruction in the last step therefore determines an approximation to the density $f'$ on the neighborhood (or patch) $U_l$. Precisely, we will implement Algorithm \ref{alg3} estimate densities on low dimensional manifolds embedded in $\mathbb{R}^d$.
\begin{algorithm}[!h]
\SetAlgoLined
\KwResult{PCA patching density estimation}
\hspace*{\algorithmicindent} \textbf{Input:} A set of test inputs $\{z'_i\}_{i=1}^{m'}\subset\mathbb{R}^d$ and a set of training inputs $X'\subset\mathbb{R}^d$\\
  \textbf{Output:} A set of density values $f'(z'_i)\in \mathbb{R}$ for $1\leq i\leq m'$ \\
\textbf{Procedure:}\\
Set the transform ($T$), either spherical ($T=T_1$) or half space ($T=T_2$)\\
Set the radius $r>0$ and the variance percentage $p$\\
\For{$i\leftarrow 1$ \KwTo $m'$}{
Set $S'=X'\cap B^d_r(z'_i)$\\
Center $S'$ so that the elements of $S'$ have mean $0$\\
Set $S\subset\mathbb{R}^n$ as the first $n$ principle component scores of $S'$, choosing $n$ to be large enough to account for $p\%$ of the variance in $S'$\\
Set $z_i\in\mathbb{R}^n$ as the transformation of $z'_i$ to principle component space\\
Reconstruct the density $h$ of the point cloud $S$ using Algorithm \ref{alg1}, with the transform set to $T$\\
Output the corresponding density value $f'(z'_i)=h(z_i)$\\
 }
 \caption{PCA patching density estimation}
 \label{alg3}
\end{algorithm}

\begin{figure}[!h]
\centering
\begin{subfigure}{0.35\textwidth}
\includegraphics[width=1\linewidth, height=4.8cm]{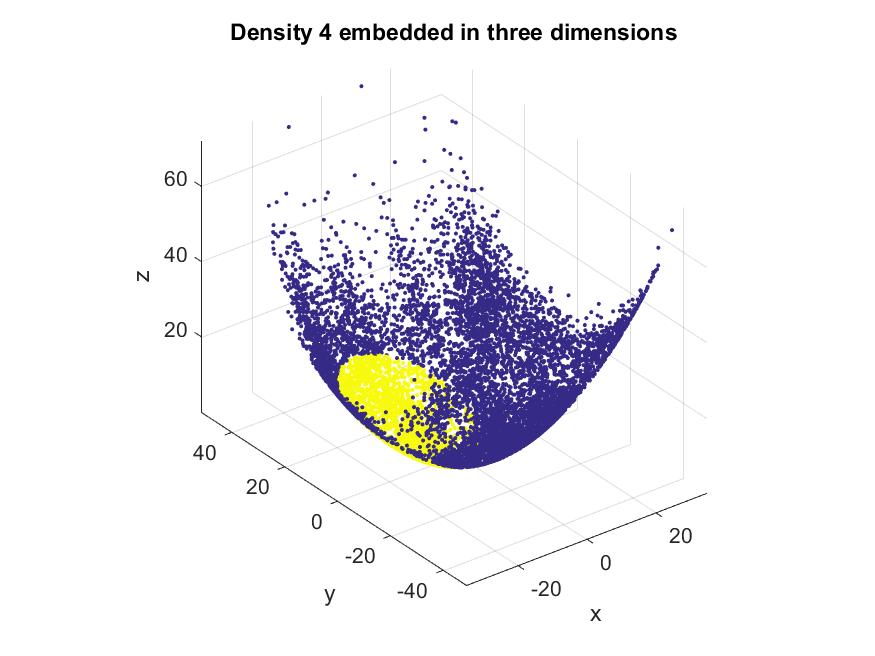} 
\end{subfigure}
\begin{subfigure}{0.35\textwidth}
\includegraphics[width=1\linewidth, height=4.8cm]{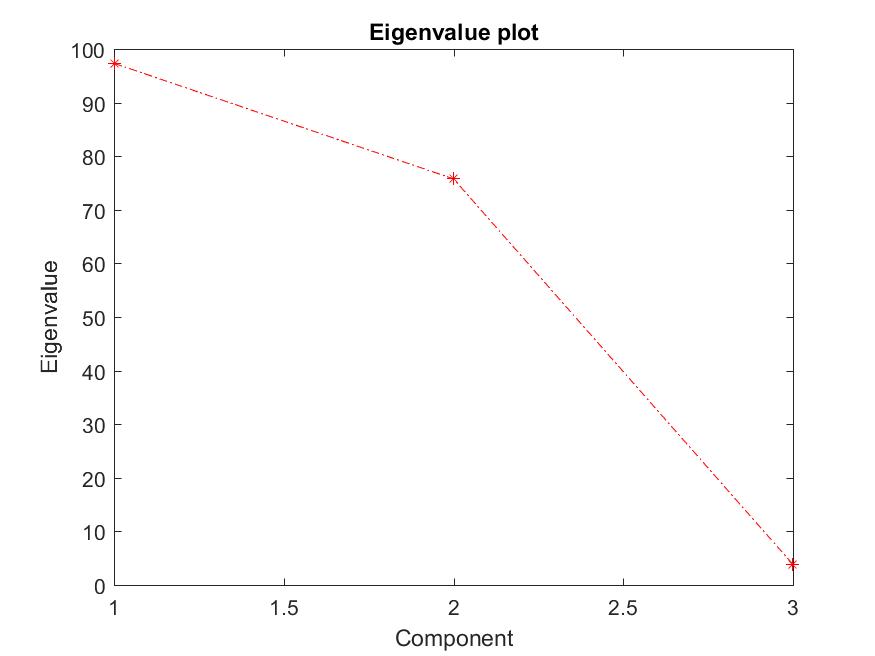}
\end{subfigure}
\begin{subfigure}{0.35\textwidth}
\includegraphics[width=1\linewidth, height=4.8cm]{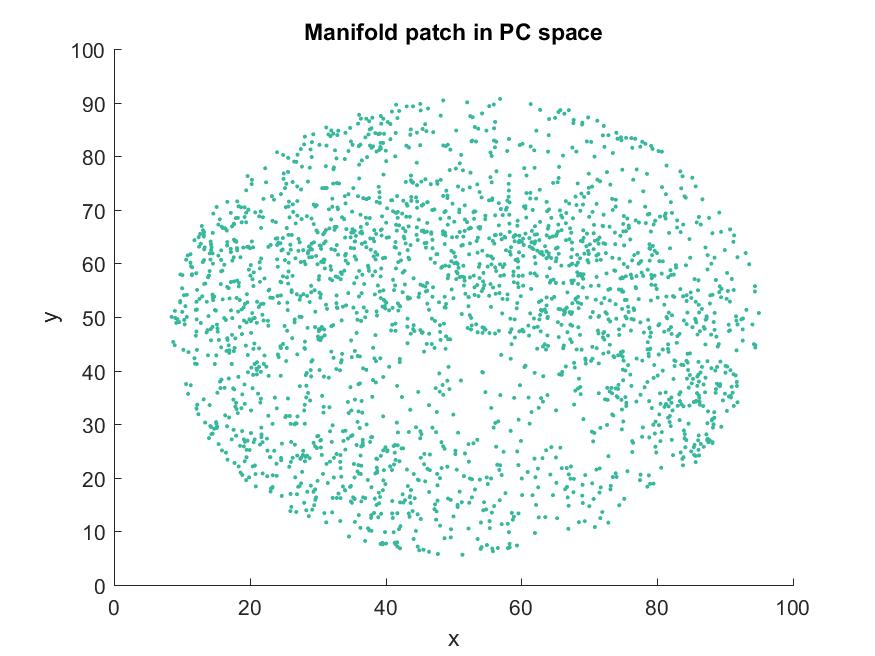} 
\end{subfigure}
\begin{subfigure}{0.35\textwidth}
\includegraphics[width=1\linewidth, height=4.8cm]{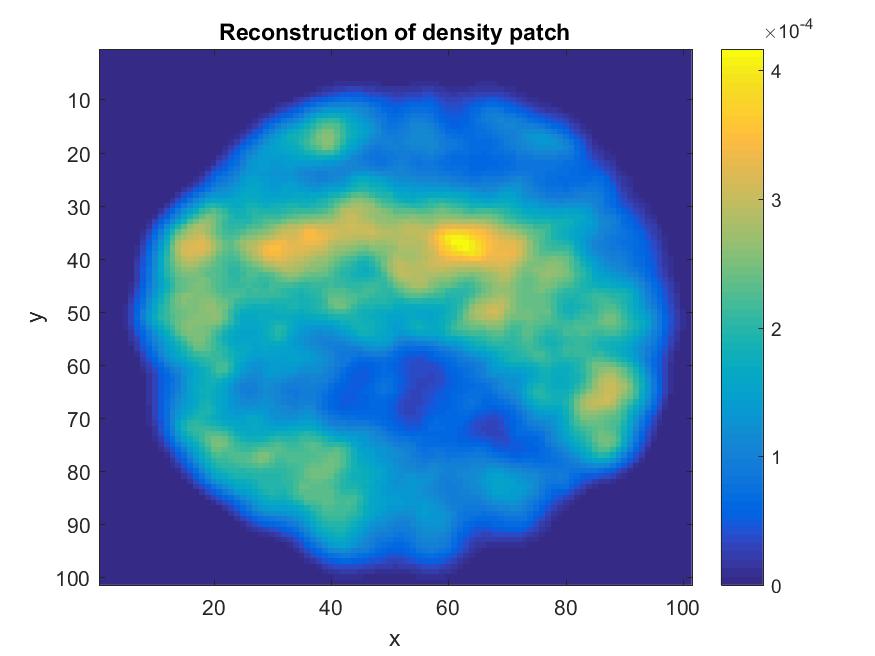}
\end{subfigure}
\caption{Density 4 embedded in three dimensions (top left) with a patch $U_l$ highlighted. The principle component eigenvalues (top right) of the patch subset. The patch density in principle component space (bottom left) and a reconstruction using spherical Radon transforms and TV (bottom right).}
\label{fig9}
\end{figure}

To demonstrate this consider figure \ref{fig9}, where we have embedded an IID sample $X$ of density 1 ($f : \mathcal{M}\to\mathbb{R}$, where $\mathcal{M}=[0,1]^2$) into three dimensions via
\begin{equation}
\phi(x,y)=(x,y,\kappa(x^2+y^2)/2),
\end{equation}
 so $f'=f \circ \phi^{-1} : \phi(\mathcal{M})\to\mathbb{R}$ is the representation of $f$ in three dimensions and $X'=g(X)$. Here $\kappa>0$ is a bound for the principle curvatures of the manifold at the origin. Note that a direct application of Algorithm \ref{alg1} here would be ill advised. Indeed if we extend $f'$ to be zero for $x\in\mathbb{R}^3\backslash\phi(\mathcal{M})$, then $R_Mf'=R_Hf'=0$ (the support of $f$ has measure zero in $\mathbb{R}^3$). In this case ($\kappa=0.05$ and $r=20$) the curvature of $\mathcal{M}$ is low on the patch $U_l=\mathcal{M}\cap B_r^3(x)$ ($x$ is located at the center of the patch) highlighted and the data dimension is recovered correctly ($n=2$), setting $p=90\%$. After transforming $S'=X' \cap B_r^d(x)$ (a density sample on $U_l$) to principle component space we can reconstruct the density on $U_l$ by applying Algorithm \ref{alg1}. See table \ref{tabEmb} for the relative errors in reconstructions of the density on $U_l$ for a range of $r$ and $\kappa$ values (keeping the $x$ and $y$ coordinates of the center of $U_l$ constant). Here we measure the relative error using equation (\ref{errmeas}) as in section \ref{method}. We see that the relative error increases with the curvature $\kappa$. This is as expected by Theorem \ref{Tspace} as the local tangent space approximations of the manifold become less accurate. We can attempt to deal with this by reducing the ball radius $r$, to counter balance the increase in the upper bound $\kappa r^2$ on the error in local tangent space approximations for increased $\kappa$. However a reduction in $r$ also serves to the lower the sample size $m$ on $U_l$ (see table \ref{tabEmb3}), and hence the reconstruction error is further amplified (by Corollary \ref{errcorr}), which is evident from the  figures in table \ref{tabEmb}. See also table \ref{tabEmb1} for the errors in the reconstruction of the density on $U_l$ using the half space Radon transform.

\section{Conclusions and further work}
Here we have presented a new approach to non-parametric density estimation inspired by ideas and theory in image reconstruction and Radon transforms. Using classical theory in geometric inverse problems we derived error estimates for our approach in section \ref{error1}. This combined the classical theory of \cite{DKW} with that of Natterer in \cite{natterer}, where we showed that the expected convergence rate to the solution was $(1/m)^{1/2(n+2)}$, where $m$ is the number of observations and $n$ is the dimension. We then went on to develop new density estimation methods in section \ref{method} from known numerical reconstruction algorithms in image reconstruction and inverse problems. The implementation of which was described in Algorithm \ref{alg1}. In particular we found that a discrete inversion of the spherical Radon transform with a total variation regularizer produced satisfactory piecewise smooth approximations to four example densities (with varying amounts of edges and smooth features)  in low dimensions ($n=2$), and offered a consistently improved performance when compared to the techniques of \cite{rad2,TV1D1} and a kernel estimate. In particular, the spherical Radon method was the most effective in recovering singularities in the density while retaining the smooth features. The half space techniques presented were also shown to offer a good performance in most cases. However we discovered that in some cases Algorithm \ref{alg1} broke down using half space transforms. This was indicated by the results of table \ref{T2}. 

We extended our methods to reconstruct densities on manifolds which are embedded in $\mathbb{R}^d$ in section \ref{manifolds}, using ideas in manifold learning and the theory of \cite{fefferman}. The implementation is described by Algorithm \ref{alg3}, which is an extention of Algorithm \ref{alg1} to local tangent spaces of manifolds. We gave an analysis of the errors in density reconstructions on manifold patches for varying levels of the principle curvature $\kappa$ (at the origin) and patch radius $r$. Our analysis was consisitent with the theory of \cite{fefferman}, and the higher curvature manifolds proved more troublesome for an accurate density estimation. To deal with this we propose to apply methods in non-linear dimensionality reduction \cite{isomap,KPCA,lle} to better ``unravel" the manifold patches of higher curvature and compare the performance to a linear approximation as is presented here.

We note that the current work is only suited to density estimation on low dimensional manifolds (mainly due to the practical restrictions of discretizing Radon transform operators in high dimensions). To deal with this, in further work we aim to develop a filtered backprojection approach for high dimensional data. As the filtering process described by an explicit inversion from hyperplane Radon transform data is ill--posed (and this is more pronounced in high dimensions), it is desired to develop a more stable filtering process that is effective in edge reconstruction, while offering little amplification in the noise. For spherical Radon transforms we aim to apply the formulae derived by John \cite[page 82]{Fritz} for expressing densities in terms of their iterated spherical means, which again will require sufficient regularization. 

Throughout this paper, we have chosen to focus on the use of empirical distribution functions to approximate spherical and half space Radon transform data in this paper. For further work we aim to use more advanced methods in density estimation (such as kernel estimators and those in \cite{TV1D,TV1D1}) to more accurately estimate the projection densities of Radon transforms, and test for an improvement in the results presented here.

\bibliographystyle{siamplain}

\clearpage
\appendix
\section[Proof of Thm]{Proof of Theorems \ref{errthm} and \ref{sphthm}}
\label{sec:proof}
Here we provide proof of Theorems \ref{errthm} and \ref{sphthm}. First some preliminary results. The next result is a statement of the Dvoretzky-Kiefer-Wolfowitz (DKW) inequality \cite{DKW} which explains the error in an approximation of a univariate density from its empirical distribution function.
\begin{theorem}[DKW theorem]
\label{DWKthm}
Let $x_1,\ldots,x_m$ be continuous, independant and identically distributed random variables in $\mathbb{R}$ with cumulative distribution function $F$. Let 
$$F_m(x)=\frac{1}{m}\sum_{i=1}^{m}I_{[x_i,\infty)}(x)$$
be the empirical distribution function for $x_1,\ldots,x_m$, where $I_S$ is the indicator function of a set $S$. Then
\begin{equation}
\label{equ15}
\mathbb{P}(\|F_m-F\|_{L^\infty(\mathbb{R})}\geq \epsilon)\leq 2e^{-2m\epsilon^2}
\end{equation}
for any $\epsilon>0$.
\end{theorem}
Next we state Bonferroni's inequality \cite[page 13]{bonferroni}.
\begin{theorem}[Bonferroni inequality]
Let $\{A_i\}_{i=1}^K$ be a finite set of $K$ events in a probability space $(\Omega,A,\mathbb{P})$, so that $\{A_i\}_{i=1}^K\subset A$. Then
\begin{equation}
\mathbb{P}\left(\cap_{i=1}^KA_i\right)\geq\sum_{i=1}^K\mathbb{P}(A_i)-(K-1).
\end{equation}
\end{theorem}
We now restate Theorem \ref{errthm} and provide proof thereafter.
\begin{theorem}
Let $x_1,\ldots,x_m$ be continuous, independant and identically distributed random variables in $\Omega^n$ with probability density function $f \in C^{\infty}_0(\Omega^n)$ and let $$g_{m}(s,\theta)=\frac{1}{m}\sum_{i\in\{x_i\cdot\theta\leq s\}}1$$ be an approximation to the half space Radon transform $g=R_Hf$. Let $K$ projection directions $\{\theta_j\}_{j=1}^K$ be uniformly spread over $S^{n-1}$. Then
\begin{equation}
\|g-g_{m}\|^2_{L^2(Z)}\leq w_{1,n-1}\frac{\log\frac{2K}{p}}{m}+\epsilon(K)
\end{equation}
with probability $1-p$ for any $0\leq p\leq 1$, where $g_{\theta}(s)=g(s,\theta)$ and $g_{m,\theta}(s)=g^{\epsilon}(s,\theta)$ and 
$$\epsilon(K)=\left|\|g-g_m\|^2_{L^2(Z)}-\sum_{j=1}^Kw_j\|g_{\theta_j}-g_{m,\theta_j}\|^2_{L^2(\mathbb{R})}\right|$$
is the error term in approximating the integral over $S^{n-1}$ by a Reimann sum, where the area elements $w_j$ are such that  $\sum_{j=1}^Kw_j=w_{1,n-1}$ and $w_{1,n-1}$ is the surface area of $S^{n-1}$.
\begin{proof}
For every $\theta\in S^{n-1}$, $g_{\theta}\in C^{\infty}_0([-1,\infty))$ is the cumulative distribution function of the one dimensional density $Rf_{\theta}\in C^{\infty}_0([-1,1])$, where $Rf_{\theta}(s)=Rf(s,\theta)$. The approximation $g_{m,\theta}=\frac{1}{m}\sum_{i\in\{x_i\cdot\theta\leq s\}}1$ is the empirical distribution function of $g_{\theta}$. After rearranging equation (\ref{equ15}) the DKW theorem gives
\begin{equation}
\mathbb{P}\left(\|g_{\theta}-g_{m,\theta}\|_{L^{\infty}(\mathbb{R})}\leq\sqrt{\frac{\log\frac{2}{p}}{2m}}\right)\geq1-p
\end{equation}
for any probability $p$ and $\theta\in S^{n-1}$. From which it follows that
\begin{equation}
\begin{split}
\|g_{\theta}-g_{m,\theta}\|^2_{L^2(\mathbb{R})}&=\int_{-1}^{1}(g_{\theta}(s)-g_{m,\theta}(s))^2\mathrm{d}s\\
 &\leq 2\|g_{\theta}-g_{m,\theta}\|^2_{L^{\infty}(\mathbb{R})}\leq\frac{\log\frac{2}{p}}{m}
\end{split}
\end{equation}
with probability $1-p$. Given $K$ projections $\{\theta_j\}_{j=1}^K$ we have
\begin{equation}
\|g_{\theta_j}-g_{m,\theta_j}\|^2_{L^2(\mathbb{R})}\leq\frac{\log\frac{2K}{p}}{m}
\end{equation}
with probability $1-p$, for all $1\leq j\leq K$ by the Bonferroni inequality. Let the $\theta_j$ be distributed uniformly over $S^{n-1}$. Then we have
\begin{equation}
\begin{split}
\|g-g_m\|^2_{L^2(Z)}&=\int_{S^{n-1}}\|g_{\theta}-g_{m,\theta}\|^2_{L^2(\mathbb{R})}\mathrm{d}\theta\\
&\leq\sum_{j=1}^Kw_j\|g_{\theta_j}-g_{m,\theta_j}\|^2_{L^2(\mathbb{R})}+\epsilon(K)\\
&\leq\frac{\log\frac{2K}{p}}{m}\sum_{j=1}^Kw_j+\epsilon(K)=w_{1,n-1}\frac{\log\frac{2K}{p}}{m}+\epsilon(K)
\end{split}
\end{equation}
with probability $1-p$, where $\epsilon(K)$ is the error in approximating the integral in the first step above by a finite sum, and the $w_j$ are small area elements of $S^{n-1}$, where $\sum_{j=1}^Kw_j=w_{1,n-1}$.
\end{proof}
\end{theorem}
We now restate Theorem \ref{sphthm} and provide proof of the result.
\begin{theorem}
\label{sphthm1}
Let $x_1,\ldots,x_m$ be continuous, independant and identically distributed random variables in $\Omega^n$ with probability density function $f \in C^{\infty}_0(\Omega^n)$ and let $$g_m(s,x)=\frac{1}{m}\sum_{i\in\{(x_i-x)^2\leq s^2\}}1$$ be an approximation to the spherical Radon transform $g=R_Mf$. Then
\begin{equation}
\frac{1}{K}\sum_{j=1}^{K}\|g_{x_j}-g_{m,x_j}\|^2_{L^2(\mathbb{R}^+)}\leq\frac{1}{m}\left(\log(K)+\log\left(\frac{2}{p}\right)\right)
\end{equation}
with probability $1-p$, for any finite set $\{x_j\}_{j=1}^K$ of circle centres. Here $g_{x_j}(s)=g(s,x_j)$ and $g_{m,x_j}(s)=g_m(s,x_j)$.
\begin{proof}
For each $x_j$, let us extent $g_{x_j}$ and $g_{m,x_j}$ to the real line by defining $g_{x_j}(s)=g_{m,x_j}(s)=0$ for $s<0$. Then $g_{m,x_j}$ is the empirical distribution function of $g_{x_j}$ and we have
\begin{equation}
\label{equ7}
\|g_{x_j}-g_{m,x_j}\|^2_{L^2(\mathbb{R})}\leq2\|g_{x_j}-g_{m,x_j}\|^2_{L^{\infty}(\mathbb{R})}\leq\frac{\log\frac{2K}{p}}{m}
\end{equation}
with probability $1-p$, for all $1\leq j\leq K$ by the DKW theorem and Bonferroni's inequality. The result follows.
\end{proof}
\end{theorem}

\section{Results on the stability of Radon tranforms}
\label{appB}
Here we provide background theory and results regarding the stability of inverse Radon transforms. First we have some definitions.
\begin{definition}
Let $f\in L^2(\mathbb{R}^n)$. Then we define the Fourier transform of $f$
\begin{equation}
\hat{f}(\xi)=(2\pi)^{-n/2}\int_{\mathbb{R}^n}f(x)e^{-i x\cdot\xi}\mathrm{d}x.
\end{equation}
\end{definition}
From \cite[page 59]{adams} we have the definitions of Sobolev spaces.
\begin{definition}
Let $\alpha\in\mathbb{R}$. Then we define the Sobolev spaces $H^{\alpha}(\mathbb{R}^n)$ of degree $\alpha$
\begin{equation}
H^{\alpha}(\mathbb{R}^n)=\left\{\text{tempered distributions}\ f : (1+|\xi|^2)^{\alpha/2}\hat{f}(\xi)\in L^2(\mathbb{R}^n)\right\}
\end{equation}
with the norm
\begin{equation}
\|f\|^2_{H^{\mathbb{R}^n}}=\int_{\mathbb{R}^n}(1+|\xi|^2)^{\alpha}|\hat{f}(\xi)|^2\mathrm{d}\xi
\end{equation}
for functions on the cylinder $Z$, $g\in H^{\alpha}(Z)$ if the norm
\begin{equation}
\|g\|_{H^{\alpha}(Z)}^2=\int_{S^{n-1}}\int_{\mathbb{R}}(1+\sigma^2)^{\alpha}|\hat{g}(\sigma,\theta)|^2\mathrm{d}s\mathrm{d}\theta
\end{equation}
exists and is finite. Here the Fourier transform is taken with respect to the variable $s$.

Let $\Omega$ be an open subset of $\mathbb{R}^n$. Then we define
\begin{equation}
H^{\alpha}(\Omega)=\{f\in H^{\alpha}(\mathbb{R}^n) : \text{supp}(f)\subseteq\bar{\Omega}\}.
\end{equation}
\end{definition}
Next from \cite[page 11]{natterer} we have
\begin{theorem}[Fourier slice theorem]
\label{FSthm}
Let $f\in C^{\infty}_0(\mathbb{R}^n)$. Then
\begin{equation}
\widehat{Rf}(\sigma,\theta)=(2\pi)^{(n-1)/2}\hat{f}(\sigma\theta),
\end{equation}
where the Fourier transform is taken with respect to the $s$ variable.
\end{theorem}
From Natterer \cite[page 42]{natterer}, we have:
\begin{theorem}
\label{Rstable}
For each $\alpha\in\mathbb{R}$ there exist positive constants $c(\alpha,n)$, $C(\alpha,n)$ such that for $f\in C^{\infty}_0(\Omega^n)$
\begin{equation}
c(\alpha,n)\|f\|_{H^{\alpha}(\Omega^n)}\leq \|Rf\|_{H^{\alpha+(n-1)/2}(Z)}\leq C(\alpha,n)\|f\|_{H^{\alpha}(\Omega^n)},
\end{equation}
\end{theorem}
Now we prove a similar result for the half space Radon transform, which is a simple consequence of the additional degree of smoothing ($+1$) applied by the integration in the $s$ variable when transforming hyperplane integral data to half space integrals.
\begin{corollary}
For each $\alpha\in\mathbb{R}$ there exist positive constants $c(\alpha,n)$, $C(\alpha,n)$ such that for $f\in C^{\infty}_0(\Omega^n)$
\begin{equation}
c(\alpha,n)\|f\|_{H^{\alpha}(\Omega^n)}\leq \|R_Hf\|_{H^{\alpha+(n+1)/2}(Z)}\leq C(\alpha,n)\|f\|_{H^{\alpha}(\Omega^n)},
\end{equation}
\begin{proof}
Since $\frac{\mathrm{d}}{\mathrm{d}s}R_Hf(s,\theta)=Rf(s,\theta)$, it is clear that $\widehat{R_Hf}(\sigma,\theta)=i\sigma\widehat{Rf}(\sigma,\theta)$, where the Fourier transform is taken with respect to the $s$ variable. From which it follows that
\begin{equation}
\begin{split}
\|R_Hf\|^2_{H^{\alpha+(n+1)/2}(Z)}&=\int_{S^{n-1}}\int_{\mathbb{R}}(1+\sigma^2)^{\alpha+(n+1)/2}|\widehat{R_Hf}(\sigma,\theta)|^2\mathrm{d}\sigma\mathrm{d}\theta\\
&=\int_{S^{n-1}}\int_{\mathbb{R}}\frac{(1+\sigma^2)^{\alpha+(n+1)/2}}{\sigma^2}|\widehat{Rf}(\sigma,\theta)|^2\mathrm{d}\sigma\mathrm{d}\theta\\
&\geq\int_{S^{n-1}}\int_{\mathbb{R}}(1+\sigma^2)^{\alpha+(n-1)/2}|\widehat{Rf}(\sigma,\theta)|^2\mathrm{d}\sigma\mathrm{d}\theta\\
&=2(2\pi)^{(n-1)/2}\int_{S^{n-1}}\int_{\mathbb{R}}(1+\sigma^2)^{\alpha+(n-1)/2}|\hat{f}(\sigma\theta)|^2\mathrm{d}\sigma\mathrm{d}\theta,\\
\end{split}
\end{equation}
where the last step follows from Theorem \ref{FSthm}. Substituting $\xi=\sigma\theta$ we get
\begin{equation}
\begin{split}
\|R_Hf\|^2_{H^{\alpha+(n+1)/2}(Z)}&=2(2\pi)^{(n-1)/2}\int_{S^{n-1}}\int_{\mathbb{R}}|\xi|^{1-n}(1+|\xi|^2)^{\alpha+(n-1)/2}|\hat{f}(\xi)|^2\mathrm{d}\xi\\
&\geq2(2\pi)^{(n-1)/2}\int_{S^{n-1}}\int_{\mathbb{R}}(1+|\xi|^2)^{\alpha}|\hat{f}(\xi)|^2\mathrm{d}\xi\\
&=2(2\pi)^{(n-1)/2}\|f\|^2_{H^{\alpha}_0(\Omega^n)},\\
\end{split}
\end{equation}
which proves the left hand inequality. The right hand inequality is left to the reader.
\end{proof}
\end{corollary}
Next we have the interpolation inequality for Sobolev spaces \cite[page 203]{natterer}.
\begin{theorem}[Interpolation inequality]
Let $f\in C^{\infty}_0(\Omega^n)$. Then
\begin{equation}
\|f\|_{H^{\lambda}}\leq\|f\|^{\frac{\beta-\lambda}{\beta-\alpha}}_{H^{\alpha}}\|f\|^{\frac{\gamma-\alpha}{\beta-\alpha}}_{H^{\beta}}
\end{equation}
for any $\alpha\leq\lambda\leq\beta$.
\end{theorem}
The next theorem gives bounds for the least squares error in a reconstruction from Radon transform data in terms of the least squares error in the Radon transform approximation. The proof of which \cite[page 94]{natterer} follows from Theorem \ref{Rstable} and the interpolation inequality.
\begin{theorem}
\label{Rerrthm}
Let $f\in C^{\infty}_0(\Omega^n)$, let $g=Rf$, and let $g^{\epsilon}$ be such that $\|g-g^{\epsilon}\|_{L^2(Z)}<\epsilon$. Then there exists a constant $c(n)$ such that
\begin{equation}
\|f-f^{\epsilon}\|_{L^2(\mathbb{R}^n)}\leq c(n)\epsilon^{1/n}\rho^{1-1/n},
\end{equation}
where $\|f\|_{H^{\frac{1}{2}}}\leq\rho$ and $f^{\epsilon}$ is such that $\|Rf^{\epsilon}-g^{\epsilon}\|_{L^2(Z)}\leq\epsilon$.
\end{theorem}
For the half space Radon transform the theorem above reads the same except for the difference in the exponent
\begin{theorem}
\label{Rerrthm1}
Let $f\in C^{\infty}_0(\Omega^n)$, let $g=R_Hf$, and let $g^{\epsilon}$ be such that $\|g-g^{\epsilon}\|_{L^2(Z)}<\epsilon$. Then there exists a constant $c(n)$ such that
\begin{equation}
\label{equ3}
\|f-f^{\epsilon}\|_{L^2(\mathbb{R}^n)}\leq c(n)\epsilon^{1/(n+2)}\rho^{1-1/(n+2)},
\end{equation}
where $\|f\|_{H^{\frac{1}{2}}}\leq\rho$ and $f^{\epsilon}$ is such that $\|R_Hf^{\epsilon}-g^{\epsilon}\|_{L^2(Z)}\leq\epsilon$.
\end{theorem}

We also have the result from Natterer \cite[page 95]{natterer}, which explains the error due to finite sampling of Radon transform data:
\begin{theorem}
\label{sample}
Let $\theta_1,\ldots,\theta_k$ and $s_1,\ldots,s_l$ be a finite set of directions and real numbers for which the pairs $(s_i,\theta_j)$ cover $S^{n-1}\times [-1,1]$ uniformly in the sense that
\begin{equation}
\sup_{-1\leq s\leq 1}\inf_{i}|s_i-s|\leq h,\ \ \ \sup_{\theta\in S^{n-1}}\inf_{j}|\theta_j-\theta|\leq h/\pi
\end{equation}
and let
\begin{equation}
d^{\mathbb{R}}(h,\rho)=\sup\left\{\|f\|_{L_2(\Omega^n)} : Rf(\theta_j,s_i)=0, \|f\|_{H^{\beta}_0(\Omega^n)}\leq\rho\right\}
\end{equation}
be the worst case error in a reconstruction of $f$ from $Rf$, where $R$ is the hyperplane Radon transform, given at the finite set of points $(s_i,\theta_j)$. Let $\beta>1/2$. Then there is a constant $c(\beta,n)$ such that
\begin{equation}
d^{\mathbb{R}}(h,\rho)\leq c(\beta,n) h^{\beta}\rho.
\end{equation}
\end{theorem}
So a uniform sampling of $Z$ is preferred to obtain the best solution using Radon transform data. Typically in tomography applications, X-ray scanners are designed to satisfy a uniform sampling rate of which there are limitations due to the hardware. 
In the current work we can sample $Rf$ as we please (although using higher $k$ and $l$ will be more computationally expensive).
\section{Additional tables and reconstruction errors}
\label{addtab}
Here we give additional tables showing the errors in reconstructions of density patches for varying levels the principle curvature $\kappa$ (at the origin) and ball radius $r$. We also give the sample sizes $m$ on the patch $U_l$ (as is represented in figure \ref{fig9}) in table \ref{tabEmb3}, for all combinations of $\kappa$ and $r$. 

\begin{figure}[!h]
\centering
\begin{subfigure}{0.42\textwidth}
\includegraphics[width=1\linewidth, height=5.5cm]{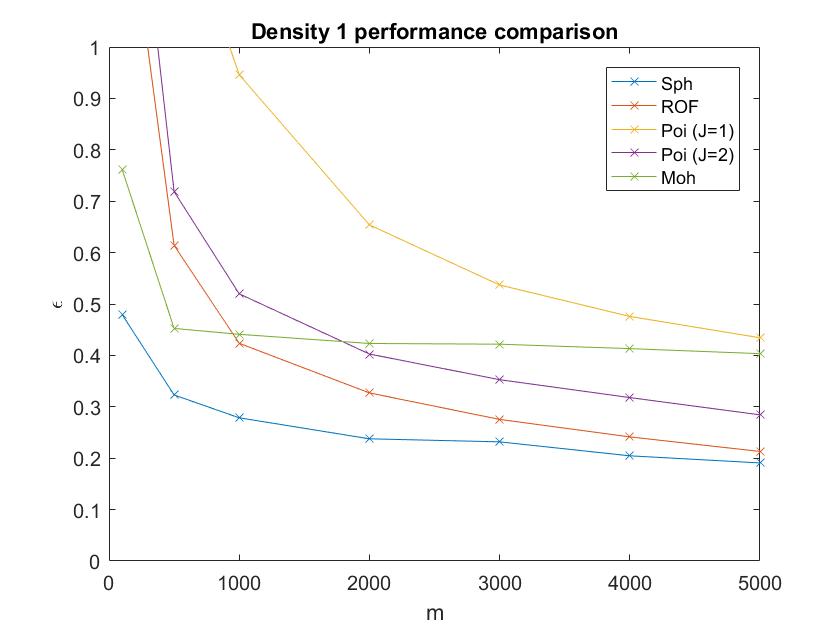} 
\end{subfigure}
\begin{subfigure}{0.42\textwidth}
\includegraphics[width=1\linewidth, height=5.5cm]{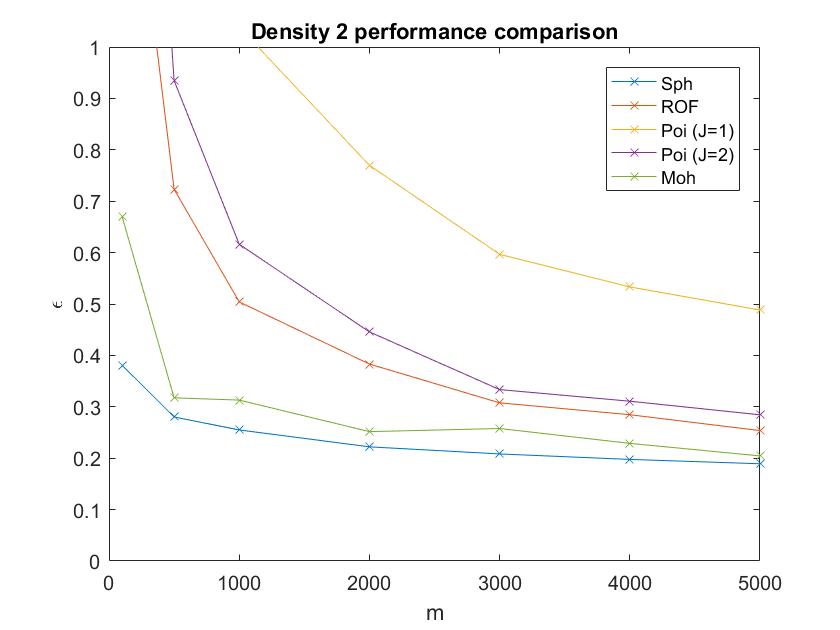}
\end{subfigure}
\begin{subfigure}{0.42\textwidth}
\includegraphics[width=1\linewidth, height=5.5cm]{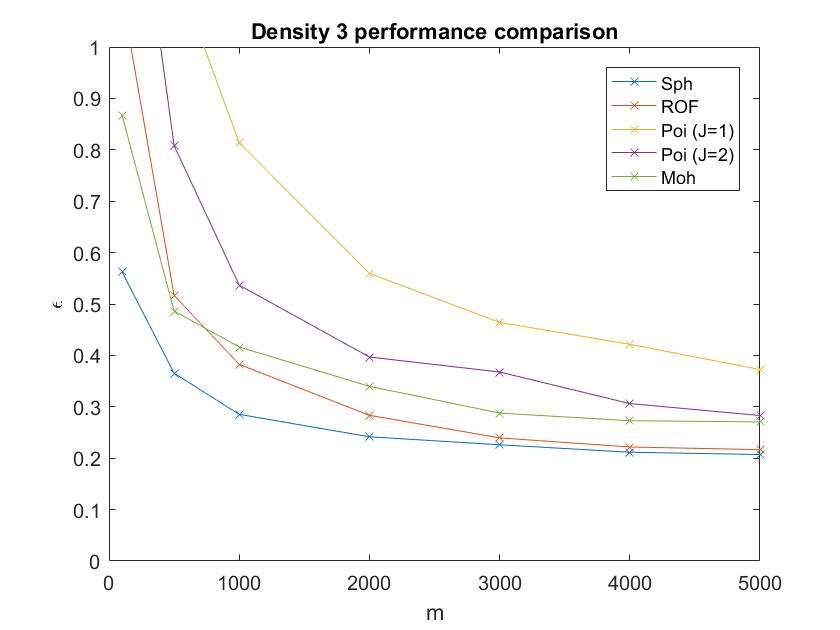} 
\end{subfigure}
\begin{subfigure}{0.42\textwidth}
\includegraphics[width=1\linewidth, height=5.5cm]{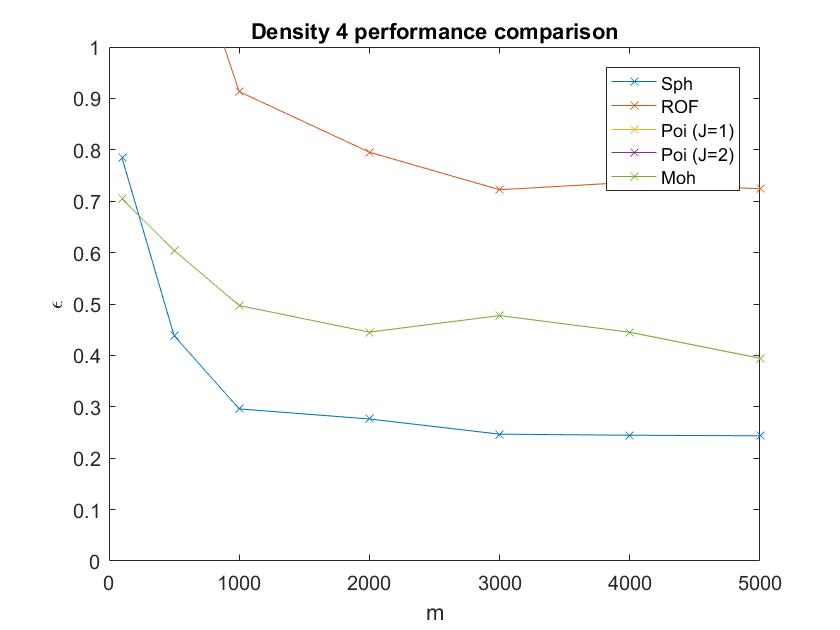}
\end{subfigure}
\caption{Comparisons of errors $\epsilon$ in density reconstructions 1--4, using the methods ROF, Poi, Sph and Moh.}
\label{figCS}
\end{figure}

In figure \ref{figCS}, we present graphs showing reconstruction errors $\epsilon$ (as in figure \ref{figC}), using TV \cite{ROF} (denoted as ROF) and Poisson \cite{Poi} (denoted as Poi) image denoising techniques. Here our noisy input image is a histogram with $100^2$ equally sized bins with centers on a 1-100 meshgrid (the histogram bins are the image pixels $p_i$). We choose the best performing (in terms of $\epsilon$) $\lambda$ from $\lambda=0.5,1,5$ (ROF) and report the results for $J=1,2$ (i.e. the $J$ such that $2^J$ divides the image dimension of 100 pixels) with LET2 (for Poi, as is shown to be most optimal in \cite{Poi}). We also report the reconstruction errors using the Split Bregman methods of Mohler {\it{et. al.}} \cite{mohler} (denoted by Moh). Here we use 10-fold CV to choose $\lambda$ (the $\lambda$ chosen is that which maximizes the sum of the log-likelihood over all folds) as in \cite{mohler}. When we ran the cross validation (using the Matlab code available in the supplementary materials of \cite{mohler}) on density 3, we experienced crashing issues upon multiple runs and on mulitple machines. Hence, for density 3, we report the minimal $\epsilon$ over the 25 $\lambda$ values considered (in \cite{mohler} they also suggest to choose from 25 $\lambda$ values). In figure \ref{figCS} we also include the reconstruction errors using Sph for comparison, and we crop the $y$ axis of the plot to only show errors between 0 and 1 (to aid the visual comparison, discounting errors greater than $100\%$). Note that for density 4, the errors using Poi were all greater than 1 (for all $m$), and hence why there are no errors curves corresponding to Poi in the bottom right hand of figure \ref{figCS}.
\begin{table}[!h]
\centering
	\begin{tabular}{| c | c | c | c | c | c | c | c |}
	\hline
		$\epsilon$ & $\kappa=0.01$ & $\kappa=0.05$ & $\kappa=0.1$     \\ \hline
$r=5$ &       $21.1\%$ &	$22.3\%$ &	$21.7\%$   \\ 
$r=10$ &       $17.4\%$ &	$21.5\%$ &	$22.5\%$   \\ 
$r=20$ &       $16.0\%$ &	$22.2\%$ &	$31.8\%$    \\ 
\hline
\end{tabular}
\caption{Relative errors in patch density reconstructions for varying levels of curvature ($\kappa$) and patch radius ($r$). The reconstructions were produced using Algorithm \ref{alg3} with the method set to half space and the variance percentage set to $p=90\%$.}
\label{tabEmb1}
\end{table}
\begin{table}[!h]
\centering
	\begin{tabular}{| c | c | c | c | c | c | c | c |}
	\hline
		$m$ & $\kappa=0.01$ & $\kappa=0.05$ & $\kappa=0.1$     \\ \hline
$r=5$ &       $196$ &	$170$ &	$118$   \\ 
$r=10$ &       $857$ &	$685$ &	$503$   \\ 
$r=20$ &       $3131$ &	$2344$ &	$1693$    \\ 
\hline
\end{tabular}
\caption{Number of samples $m$ on the patch $U_l$ for varying levels of curvature ($\kappa$) and patch radius ($r$).}
\label{tabEmb3}
\end{table}
\clearpage
\section{Comparison of spherical and half space reconstructions}
\label{appC}
\begin{figure}[!h]
\centering
\begin{subfigure}{0.42\textwidth}
\includegraphics[width=1\linewidth, height=5.5cm]{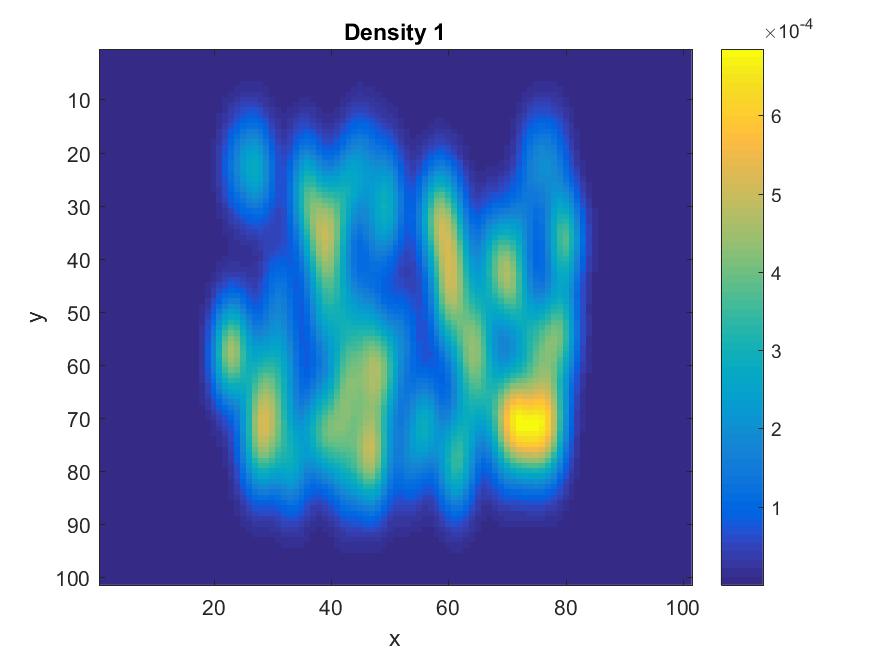} 
\end{subfigure}
\begin{subfigure}{0.42\textwidth}
\includegraphics[width=1\linewidth, height=5.5cm]{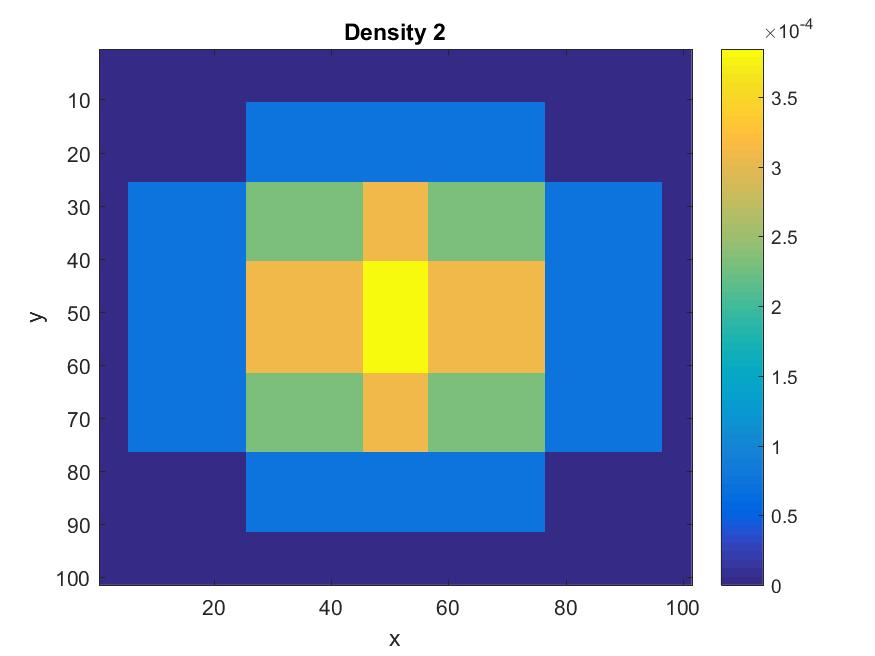}
\end{subfigure}
\begin{subfigure}{0.42\textwidth}
\includegraphics[width=1\linewidth, height=5.5cm]{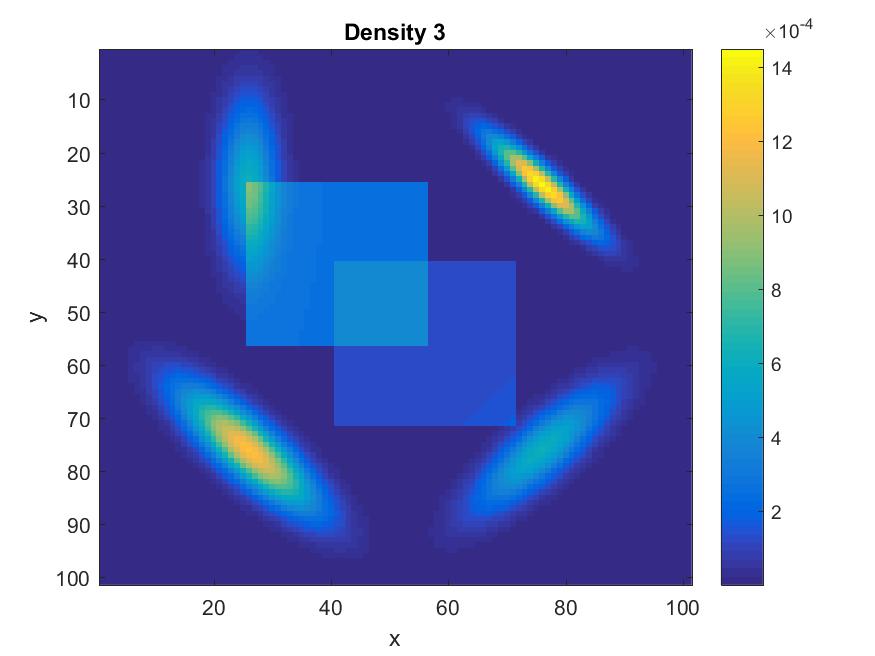} 
\end{subfigure}
\begin{subfigure}{0.42\textwidth}
\includegraphics[width=1\linewidth, height=5.5cm]{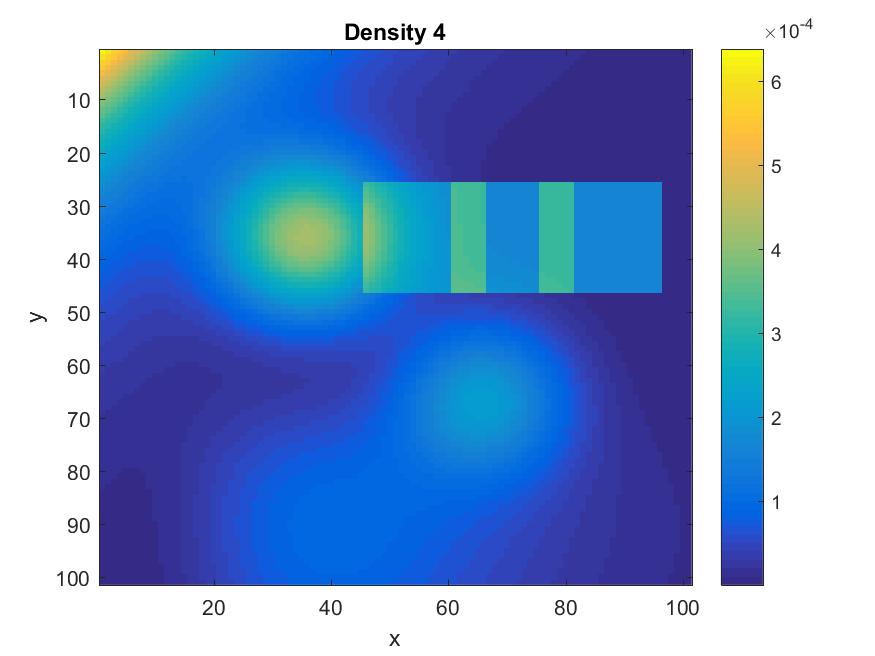}
\end{subfigure}
\caption{Ground truth density functions.  A Gaussian mixture with 100 Gaussians (density 1, top left), 5 overlapping uniform densities (density 2, top right), a Gaussian--uniform mixture density (density 3, bottom left) and a general mixture density (a sum of exponential, Gaussian, gamma and uniform density functions (density 4, bottom right)).}
\label{fig3.3}
\end{figure}

\begin{figure}
\begin{subfigure}{0.49\textwidth}
\includegraphics[width=0.9\linewidth, height=6cm]{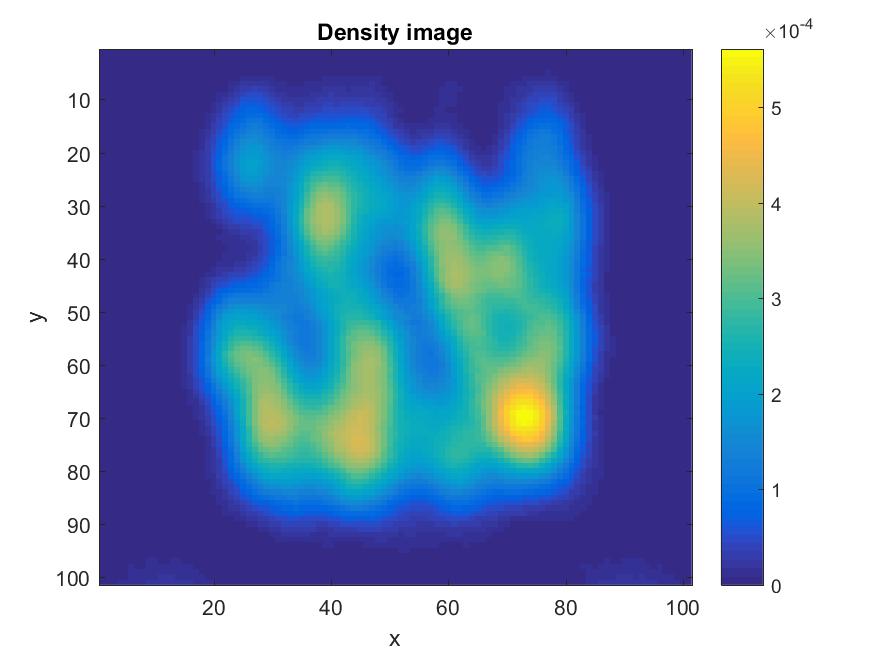} 
\end{subfigure}
\begin{subfigure}{0.49\textwidth}
\includegraphics[width=0.9\linewidth, height=6cm]{Sph1.jpg}
\end{subfigure}
\caption{Gaussian mixture reconstructions using half space (left) and spherical (right) Radon transforms, with TV.}
\label{figC1}
\end{figure}

\begin{figure}
\begin{subfigure}{0.49\textwidth}
\includegraphics[width=0.9\linewidth, height=6cm]{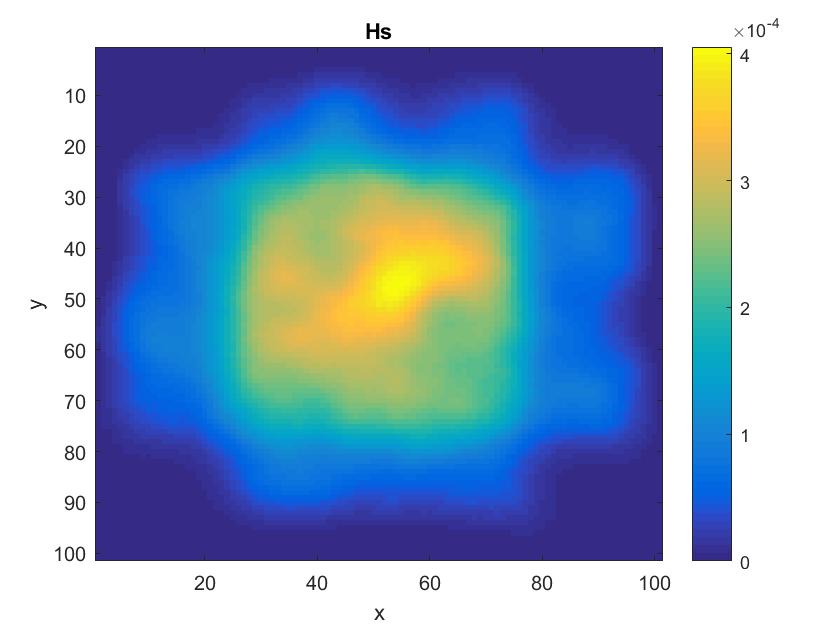} 
\end{subfigure}
\begin{subfigure}{0.49\textwidth}
\includegraphics[width=0.9\linewidth, height=6cm]{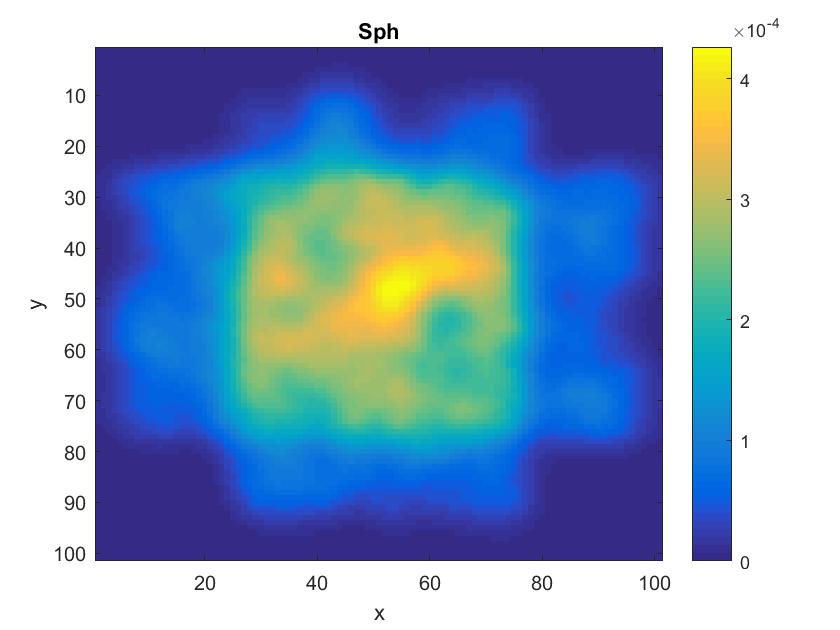}
\end{subfigure}
\caption{Density 2 reconstructions using half space (left) and spherical (right) Radon transforms, with TV.}
\label{figC2}
\end{figure}

\begin{figure}[!h]
\begin{subfigure}{0.49\textwidth}
\includegraphics[width=0.9\linewidth, height=6cm]{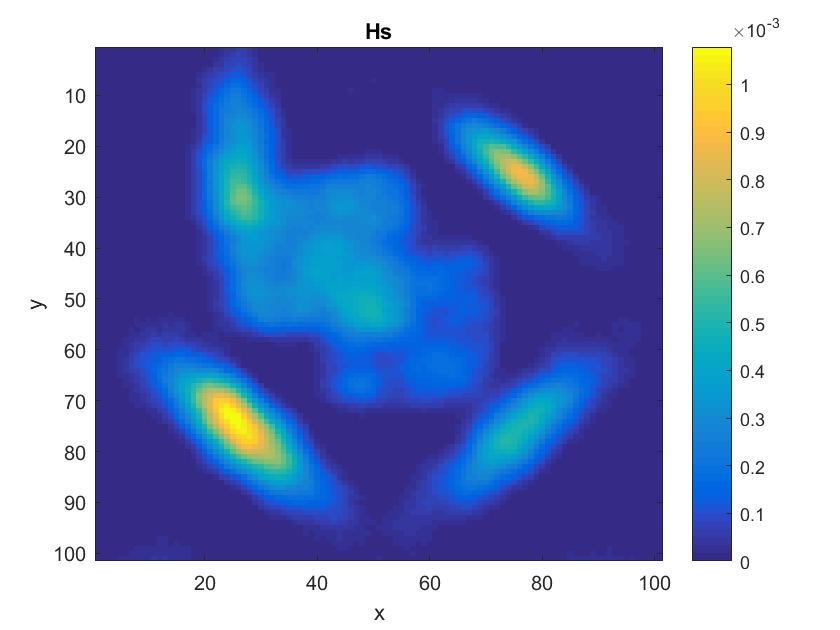} 
\end{subfigure}
\begin{subfigure}{0.49\textwidth}
\includegraphics[width=0.9\linewidth, height=6cm]{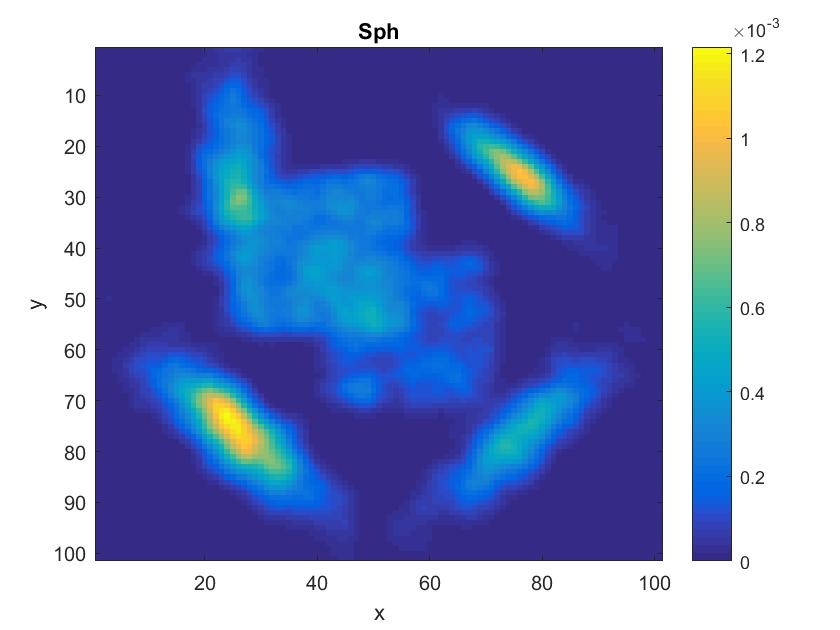}
\end{subfigure}
\caption{Density 3 reconstructions using half space (left) and spherical (right) Radon transforms, with TV.}
\label{figC3}
\end{figure}

\begin{figure}[!h]
\begin{subfigure}{0.49\textwidth}
\includegraphics[width=0.9\linewidth, height=6cm]{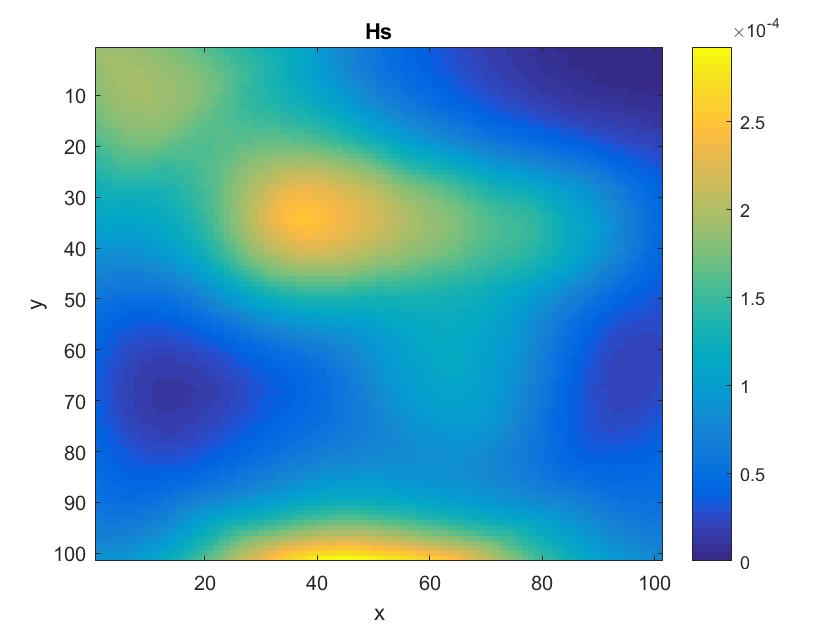} 
\end{subfigure}
\begin{subfigure}{0.49\textwidth}
\includegraphics[width=0.9\linewidth, height=6cm]{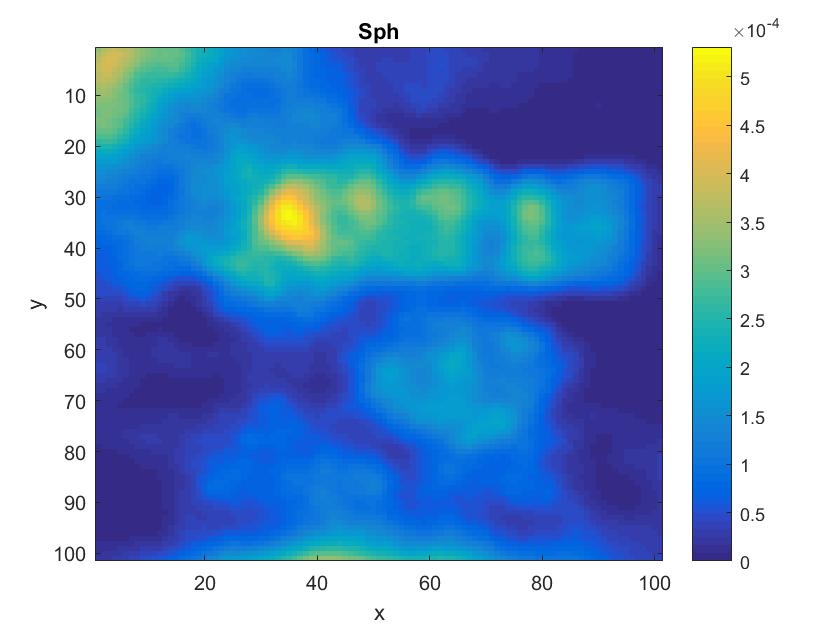}
\end{subfigure}
\caption{Density 4 reconstructions using half space (left) and spherical (right) Radon transforms, with TV.}
\label{figC4}
\end{figure}

Here we give a side by side comparison of reconstructions of densities 1--4 using spherical and half space Radon transforms, with $m=2000$ draws. Figures \ref{figC1}--\ref{figC4} give a comparison of spherical and half space results using a TV and GCV method as in section \ref{method}. We also present larger views of the ground truth densities 1--4 in figure \ref{fig3.3}, with descriptions.

When using a TV regulariser we notice significantly improved results using spherical Radon transforms for densities 1--3 (although less significant in the case of density 2), and there are significant artefacts introduced in the reconstruction of density 4 with the half space Radon transform. The spherical and half space Radon transforms are both stably invertible in some Sobolev space, so we would expect to see similar results, but in the implementation the errors in the reconstructions of densities 1--4 when using TV are improved with spherical transforms. This could be due to the over determined nature of spherical data (which implies a more stable inversion), or possibly due to the differences in the geometry of spheres and hyperplanes (e.g. one is a closed surface and the other is open). 
\clearpage


\end{document}